\documentclass[11p,leqno]{article}
\usepackage{amsmath}
\usepackage{amsfonts}
\usepackage{amssymb}
\usepackage{graphicx}
\usepackage{color}
\newtheorem{thm}{Theorem}[section]
\newtheorem{cor}[thm]{Corollary}
\newtheorem{lem}[thm]{Lemma}
\newtheorem{prop}[thm]{Proposition}

\numberwithin{equation}{section}

\newcommand{\dx}{\,{\rm d}x}
\newcommand{\dy}{\,{\rm d}y}

\newcommand{\ds}{\,{\rm d}s}
\newcommand{\dt}{\,{\rm d}t}
\newcommand{\rd}{\,{\rm d}}
\newcommand{\vol}{{\rm Vol}}

\DeclareMathOperator*{\appros}{\approx}
\DeclareMathOperator*{\asint}{\sim}
\def\LL{\mathrm{L}} %per gli spazi L^p
\newcommand{\RR}{\mathbb{R}}
\newcommand{\NN}{\mathbb{N}}

\newcommand{\ve}{\varepsilon}
\def\ee{\mathrm{e}} %per l'esponenziale
\def\qed{\unskip\kern 6pt \penalty 500
\hfuzz=3pt \raise -2pt\hbox{\vrule \vbox to8pt{\hrule width 6pt
\vfill\hrule}\vrule}\par}
\def\beq{\begin{equation}}
\def\eeq{\end{equation}}
\def\proof{\noindent {\sl Proof.~}}
\definecolor{darkblue}{rgb}{0.05, .05, .65}
\definecolor{darkgreen}{rgb}{0.05, .70, .05}
\definecolor{darkred}{rgb}{0.8,0,0}

\begin{document}
\title{\bf Special fast diffusion with slow asymptotics. \\
Entropy method and flow on a Riemannian manifold}
\author{Matteo Bonforte$^{\rm 1,3}$, Gabriele Grillo$^{\rm 2,4}$,
\ Juan Luis V\'azquez$^{\rm 1,5}$}
\date{} %%  this cancels date in article format
\maketitle

\begin{abstract}
We consider  the asymptotic behaviour  of positive solutions $u(t,x)$ of the fast diffusion
equation $u_t=\Delta \left(u^{m}/m\right)=\mbox{\rm div}\,(u^{m-1}\nabla u)$ posed for $x\in\RR^d$, $t>0$,
with a precise value for the exponent $m=(d-4)/(d-2)$. The space dimension is $d\ge 3$ so that $m<1$, and even $m=-1$ for $d=3$.
This case had been left open in the general study
\cite{BBDGV} since it requires quite different functional analytic
methods, due in particular to the absence of a spectral gap for the operator generating the
linearized evolution.

The linearization of this flow is interpreted here as the heat flow
of the Laplace-Beltrami operator of a suitable Riemannian Manifold
$(\RR^d,{\bf g})$, with a metric ${\bf g}$ which is conformal to the
standard $\RR^d$ metric. Studying the pointwise heat kernel
behaviour allows to prove {suitable Gagliardo-Nirenberg}
inequalities associated to the generator. Such inequalities in turn
allow to study the nonlinear evolution as well, and to determine its
asymptotics, which is identical to the one satisfied by the
linearization. In terms of the rescaled  representation, which is a nonlinear Fokker--Planck equation, the convergence rate turns out to be polynomial in time. This result is in contrast with
the known exponential decay of such representation for all other
values of $m$.

\vskip 1.5cm

\end{abstract}
\
\noindent (1) Departamento de Matem\'{a}ticas, Universidad
Aut\'{o}noma de Madrid, Campus de Cantoblanco, 28049 Madrid, Spain

\noindent (2) Dipartimento di Matematica, Politecnico di Torino,
corso Duca degli Abruzzi 24, 10129 Torino, Italy
%\noindent (3) Ceremade, Universit\'e Paris Dauphine, Place de Lattre
%de Tassigny, F-75775 Paris C\'edex 16, France

\noindent (3) e-mail address:~matteo.bonforte@uam.es

\noindent (4) e-mail address:~gabriele.grillo@polito.it

\noindent (5) e-mail address:~juanluis.vazquez@uam.es
%\date{} %%  this cancels date in article format
\maketitle
%%%%%%%%%%%%%%%%%%%%%%%%%%%%%%%%%%%%%%%%%%%%%%%%%%%%%%%%%%%%%%%%
%\newpage
%\tableofcontents
\newpage

%%%%%%%%%%%%%%%%%%%%%%%%%%%%%%%%%%%%%%%%%%%%%%%%%%%%%%%%%%%%%%%%%%%%%%

\section{Introduction} \label{sect.intro}

In this paper we shall describe the asymptotic behaviour (as $t\to\infty$) of a class
of solutions $u(t,x)\ge 0$ of the  fast diffusion equation (FDE)
\begin{equation}\label{1.1}
\partial_t u=\Delta\left(u^{m}/m\right)=\nabla \cdot(u^{m-1}\nabla u),
\qquad m<1,
\end{equation}
 posed\footnote{There is no restriction $m>0$. The last expression represents a parabolic equation whenever $u>0$
even if $m\le 0$. For $m=0$ the first expression must be replaced by
$\Delta\log(u)$.} for $t>0$ in the whole space, $x\in \RR^d$, in
dimensions $d\ge 3$, and taking  initial data
\begin{equation}
 u(0,x)=u_0(x)>0,
 \end{equation}
where $u_0$ belongs to a class to be made precise below, in
particular $u_0$ is bounded and decays at infinity like
$c\,|x|^{2/(1-m)}$ with lower order terms. Actually, since $m<(d-2)/d$ it is well-known
that for initial data of the above form  the weak solution exists
and is unique for small times, and then extinguishes completely after
a finite time $T=T(m,d,u_0)$, \cite{VazSmooth}. We are interested in
the behaviour of such solutions near extinction, as $t\nearrow T$. A
detailed analysis of this question has been performed in a recent paper
\cite{BBDGV} for general $m<1$ (even when $m\le 0$), but rates of
convergence could not be obtained for a special value of the
diffusion exponent $m$, precisely\footnote{In space dimension  $d=4$ we have $m_*=0$,
logarithmic diffusion. For $d=3$ we deal with $m_*<0$, a very
singular case that was only briefly exposed in \cite{BBDGV}.} for
$m_*=(d-4)/(d-2)$. We refer to that paper for further references to
the abundant literature on the topics of entropy methods, rescaling
and rates of convergence for this type of nonlinear diffusion
equations, cf. also \cite{MR1853037}, \cite{MR1986060},
\cite{Daskalopoulos-Sesum2006}, \cite{MR1940370}, \cite{DenzMcCann},
\cite{HP}, \cite{MR1974458}, \cite{Otto}, \cite{VazAs}.

The present paper is devoted to settle the asymptotic behaviour in
the special case $m=m_*$. We shall see that it falls out of the scope
of asymptotic theory developed in the paper \cite{BBDGV} for the
rest of the values $m<1$, both in the type of techniques and in the
type of results. The clue to finding the stabilization rates of the
rescaled orbits towards their equilibrium states in this special
case relies on

(i) Realizing that a suitable linearization of the rescaled flow
can be viewed as plain heat flow in a suitable Riemannian manifold.
This allows us to use the very detailed theory that has been
developed for studying (the long-time behaviour of) such flows, see
\cite{LY, D};

(ii) Performing a study of nonlinear stability based on an interesting
modification of the entropy methods of \cite{BBDGV}.

\noindent  The paper gives precise statements and proof of these
assertions. It is organized as follows: in the next section we shall
review the needed facts about the asymptotic behaviour of our
problem in the more general setting of variable $m\in \RR$. We also introduce
 the family of entropies that allow to
prove the plain stabilization result of \cite{BBDGV}, as well as the
linearization method that allows to obtain rates of convergence when
$m\ne m_*$, when used  in combination with the limit of the previous
entropies. The failure of this approach in the special case $m_*$ is
identified in  \cite{BBDGV} as the lack of a suitable {\sl spectral gap}
\ in the operator analysis of the linearized problem.

We then focus on $m=m_*$ and address such an essential difficulty.
The convergence results are carefully stated in Section \ref{ssec.statement}.
We start the new work in Section \ref{sect-lin} by a detailed
analysis of the linearized equation, identified as a heat flow on a
cigar-like Riemannian manifold. This is followed by the results on
linearized stability. Section \ref{sect.nlem} gathers all the results
needed in the comparison of linear and nonlinear entropies.
The proof of nonlinear stability is given in Section \ref{sect.nlem2}.
In Section \ref{m.neq.mstar} we revise the convergence
for the case $m\neq m_*$ and show that our method provides a shorter proof and also a small
improvement with respect to \cite{BBDGV}.

The main difference in the asymptotic results is that convergence to
a selfsimilar profile takes place with a rate of
approach that differs in a marked way from the power
rate of all the cases $m\ne m_*$. The convergence is
most clearly visualized below in the rescaled representation,
a nonlinear Fokker--Planck equation, where it takes the form of stabilization towards equilibrium
with a polynomial rate of approach in terms of the new time
variable $s$. Specifically, the study is made in terms of rescaled
variable
\begin{equation}
v(s ,y)=(T-t)^{-d\beta}u(t,x), \quad y=ax (T-t)^{\beta}, \quad s
= \gamma\log(T/(T-t)),
\end{equation}
where $T>0$ is the extinction time and the constants $\beta,\gamma$
and $a$ are precisely defined in Section \ref{sect.prel}\footnote{The exponent $\beta$ is essential, whereas
the values of $a, \gamma>0$ are just convenient.}. This
rescaled variable satisfies the nonlinear Fokker-Planck equation,
see \eqref{eq.v1} or \eqref{eq.v},   which is better suited for the asymptotic
analysis. The stationary profile for the latter version of the equation is given by the simple
expression
\begin{equation}
V_D(y)=1/(D+|y|^2)^{(d-2)/2}, \qquad D>0,
\end{equation}
for a suitable constant $D$ determined by the initial data. This simple expression is handy
since $V_D$ and powers of it will appear as weights in some functional inequalities that are
essential in our study. In terms of the new logarithmic time $s$ (that goes to
infinity as $t\to T$), the long time behaviour of the rescaled flow takes the form of
stabilization towards the  profile $V_D$ with a power rate of
convergence:
 \begin{equation}
 \|v(s,y)-V_D(y)\|_{L^\infty(\RR^d)}= O(s^{-1/4}) \quad \mbox{as \ } s\to
 +\infty.
 \end{equation}
This rate replaces the exponential decay formulas with respect to
$s$ of the cases $m\ne m_*$,  that have been obtained
in \cite{BBDGV}.   This polynomial rate in $s$ is
slower than the exponential rate in $s$ that obtains in all other
cases $m<1$, $m\ne m_*$. Summing up, we are in a case of what is
called {\sl slow asymptotics}, or {\it critical slowing down}, in mechanical
systems and statistical mechanic, and such cases
need as a rule special analytical methods.

The needed assumption on the initial data is that $v(0,y)$ be a small
perturbation of $V_D(y)$ in a sense made precise by assumptions
(H1') and (H2') below. Let us  stress that some kind of
similar assumption on the data is needed to obtain the asymptotic
result. Actually, for data that decay at infinity with a slower rate than
$O(|y|^{2/(1-m)})$ (i.e., with a smaller power) solutions do not
even extinguish in finite time. On the other hand, for data that
decrease with a larger power, the behaviour near the extinction time follows a
completely different pattern that is described in the monograph \cite{VazSmooth}.

We complete this introduction with some comments on related topics.
Let us first recall that there are two other
known instances of interpretation of fast diffusions as
geometrical flows. The first case is the evolution Yamabe flow, i.\,e., the fast diffusion with
$m=(d-2)/(d+2)$, $d\ge 3$. It describes how a conformal
Riemannian metric evolves by scalar curvature; in that case $u$ is
interpreted as the conformal factor of the metric raised to the
power $(d+2)/4$. An asymptotic study of this problem is made by Del Pino and S\'aez  in
\cite{DPS} with exponential convergence to a separate variable solution, and the results
are extended in \cite{VazSmooth}. In the second case we
deal with Ricci flow in dimension $d=2$, as proposed by
Hamilton \cite{Ham}, and then $m=0$ (logarithmic diffusion). The asymptotic
behaviour in that case is rather complex, cf. \cite{DH} or the monograph \cite{VazSmooth}. Both models happen in a different context, since they consist of
interpreting the variable $u$ in the FDE as the evolving conformal
factor of a conformal representation, while here we consider a
heat flow on a fixed manifold as the linearization limit of a
nonlinear fast diffusion flow. They have in common the property of
extinction in finite time.

Finally, we mention that a number of formulas and ideas used in the theory
of  Ricci flows bear a close similarity with developments in linear and
nonlinear  diffusion theory. Thus, the use of entropies is prominent in
Perelman's study of  the Ricci flow, \cite{Perel}, where he introduces his
 functionals $\cal F$ and $\cal W$  which are extensions of the Einstein-Hilbert functional. He then writes  the gradient flow for the functionals  as a system of equations for the evolving  metric $g_{ij}$ and a scalar
function $f$, which satisfies a backward heat  equation.  Strong connections exist with studies of entropies for heat equations on a static  manifold, see for example Ni \cite{Ni1} and also the general
 references \cite{ChowK, ChowLN, Muller}. In a recent paper \cite{LNVV}
Lu, Ni,  Villani and one of the authors investigate Harnack inequalities
and entropies for  porous medium and fast diffusion equations on static
manifolds that are closely  related  to Yau, Hamilton and Perelman's work,
 and on the other hand are close to  the  subject of this paper.
The whole topic calls for further understanding.

%%%%%%%%%%%%%%%%%%%%%%%%%%%%%%%%%%%%%%%%%%%%%%%%%%%%%%%%%%%%%%%%
\subsection*{List of notations}

\par\noindent $D_0,D_1,D_*$: the constants involved in Assumptions
(H1), (H1$^\prime$), (H2),(H2$^\prime$). See Section \ref{s2.2}.
\par\noindent $f$, $\tilde f$: the functions involved in Assumptions
(H2), (H2$^\prime$). See Section \ref{s2.2}.
\par\noindent ${\mathcal F}(w)$: the relative entropy. See Formula
\eqref{Entropy.Quotients}.
\par\noindent$F$: the linearized relative entropy. See Proposition
\ref{4.13} and Lemma \ref{Lem.Bounds.RE}. The argument of $F$ can be
both $g$ and $w$ (see below for the meaning of the latter
quantities).
\par\noindent$g$: the (weighted) linearization of $w-1$. See
Formulas \eqref{2.14} and \eqref{lin.g}.
\par\noindent ${\bf g}_\alpha$: the metric describing the geometric
interpretation of the linearized operator. See Formula
\eqref{metric}.
\par\noindent ${\mathcal I}(w)$: the relative Fisher information. See
Formula \eqref{Fisher.Quotients}.
\par\noindent $I_m$: the linearized Fisher information, see
\eqref{form}. The index $m$ is dropped in Section 5 for brevity.
\par\noindent$K(t,x,y)$: the heat kernel of the Laplace--Beltrami
associated to ${\bf g}_\alpha$. See Section \ref{linear}.
\par\noindent$L_m$: the linearized generator. See Formula \eqref{op}.
\par\noindent$\mu_*$: the weighted measure $\rd\mu_*=V_D^{2-m_*}\rd x$.
See just before Section
\ref{funct}. The L$^p$ norms in Section 4 are taken w.r.t. $\mu_*$.
\par\noindent $T$: the extinction time of the
Barenblatt solutions and of the solutions considered. See Formula
\eqref{2.3} and Assumption (H1).
\par\noindent$u(x,t)$: the solution to the fast diffusion equation.
See Formula \eqref{1.1}.
\par\noindent$U_D(t,x)$: the Barenblatt solutions for $m>m_c$. See
Formulas \eqref{2.1} and \eqref{fc}.
\par\noindent$U_{D,T}(t,x)$: the pseudo--Barenblatt solutions for
$m<m_c$. See Formula \eqref{2.3}. \par\noindent$v(y,s)$: the
rescaled solution of the nonlinear Fokker--Planck equation. See
Formula \eqref{eq.v}.
\par\noindent$V_D(y)$: the Barenblatt profiles in rescaled
variables. See Formula \eqref{2.7}. $V_*(y):=V_{D_*}(y)$ is defined
in Section \ref{entropies}.
\par\noindent$w$: the ratio $v/V_{D_*}$. See Formula \eqref{eq.w} for the
equation satisfied by $w$.
\par\noindent$W_0,W_1$: lower and upper bounds for $w$. See Section
\ref{entropies}.

\medskip

\noindent The notation $\|\cdot\|_p$ denotes in principle the
standard norm in ${\rm L}^p({\mathbb R}^d)$, but starting at the end
of Subsection \ref{sect-lin}.1 we will use weighted spaces and it
will indicate ${\rm L}^p({\mathbb R}^d,\rd \mu)$ with a weight $\mu$
related to the Barenblatt solutions. The context will always make it
clear.

%%%%%%%%%%%%%%%%%%%%%%%%%%%%%%%%%%%%%%%%%%%%%%%%%%%%%%%%%%%%%%%%
\section{Preliminaries: rescaling, stabilization and entropy}
\label{sect.prel}

The fast diffusion equation with $0<m<1$ has attracted the attention
of researchers  in recent times, once the theory of the
corresponding slow diffusion case $m>1$ came to be well known. In
the latter case the long-time behaviour of all solutions with nonnegative
and $L^1$ data $u_0$ is given by a one-parameter family of explicit
self-similar solutions of the form
\begin{equation}\label{2.1}
U_D(t,x)=t^{-\alpha}B_D(xt^{-\beta}),
\end{equation}
with $\beta=1/(2+d(m-1))$, $\alpha=d\beta$ and profile
$B_D=(D-k|\xi|^2)_+^{1/(m-1)}$ with a free constant $D>0$,  a fixed
constant $k=\beta(m-1)/2$, and putting $\xi=xt^{-\beta}$. These
solutions, usually called Barenblatt solutions, replace the Gaussian
profiles found in the long time behaviour of the classical heat
equation, which is the case $m=1$. See the precise asymptotic result
in \cite[Chapter 18]{BookPME}.

When going over to the fast diffusion equation,  the situation has
been well understood in a first range of exponents $1>m>m_c=(d-2)/d$
(the `good' fast diffusion range); indeed,  solutions of the above
initial value problem exist and are unique, they are positive and smooth for
every choice of the initial data in $\LL^1_{\rm loc}(\RR^d)$, and even in
more general cases, cf. \cite{ChassVaz}. In particular, Barenblatt
solutions still exist, they have the same selfsimilar form though
the profile looks a bit different
\begin{equation}\label{fc}
B_D(\xi)=(D+k|\xi|^2)^{-1/(1-m)}
\end{equation}
now with $k=\beta(1-m)/2$. This is a positive function everywhere
in $\RR^ d$ and decays at infinity like $O(|\xi|^{-2/(1-m)})$, so
that $B_D\in L^1(\RR^d)$ if $m>m_c$. The Barenblatt solutions still represent
the asymptotic behaviour of all solutions with nonnegative and
$L^1$ data $u_0$, with even better convergence result in relative error, cf.
\cite{BV, VazAs}. Factors like $B_D$ will appear in the sequel as
weights in functional inequalities and measure spaces. We shall use
below a proper scaling to get rid of the inessential constant $k$.

However, such a simple theory breaks down for $m<m_c$, even if $m>0$
(which is possible if $d\ge 3$), due in particular to the phenomenon
of extinction in finite time, cf. \cite{VazSmooth}. In particular,
our model solutions cannot be continued in the same form because the
similarity exponents $\alpha$ and $\beta$ go to infinity as $m$ goes
down to $m_c$. But for $m<m_c$ a related family of extinction
solutions is found of the backward self-similar form
\begin{equation}\label{2.3}
U_{D,T}(t,x)=(T-t)^{\alpha}B_D(x(T-t)^{\beta}),
\end{equation}
with $\beta=1/(d(1-m)-2)>0$ and $\alpha=d\beta>0$ (just minus the
formulas used  before). Here, $T$ and $C$ are arbitrary positive
constants and  $B_D$ is given just as in the case $m_c<m<1$. It is
to be noted that $B_D$ is no more an integrable function in $\RR^d$,
so we are completely away from the functional setting we started
from. These new solutions are sometimes called pseudo-Barenblatt
solutions to distinguish them from the original Barenblatt family.

\subsection{Rescaled flow equation}

 Actually, these
solutions do not possess the strong attractivity properties of
their relatives for $m>m_c$. In order to investigate their partial
attractivity (more precisely, their rescaled stability), we have
studied in the paper \cite{BBDGV} the extinction behaviour of
solutions with initial data close to a pseudo-Barenblatt solution.
This is the situation in short terms: we can show that after a
rescaling step we obtain the nonlinear Fokker-Planck equation
\begin{equation}\label{eq.v1}
\partial_s v= \frac{a^2}{\gamma}\nabla_y (v^{m-1}\nabla_y v)+
\frac{\beta}{\gamma}\nabla_y\cdot (y v)
\end{equation}
in terms of the rescaled variable $v(s ,y)$ defined as
\begin{equation}
v(s ,y)=(T-t)^{-d\beta}u(t,x), \quad y=ax (T-t)^{\beta}, \quad s
= \gamma\log(T/(T-t)).
\end{equation}
Here $T=T(u_0)$ is the extinction time of the solution, $\beta=(d(1-m)-2)^{-1}$ and
we will choose the free constants $a,\gamma>0$
 to be $a^2=\gamma=(1-m)\beta/2$. Note that $s (0)=0$ and $s (t)\to\infty$ as $t\to T$. This means
that whenever we use as $T$ the actual extinction time of the
solution $u$, then $v$ is globally defined, for $y\in \RR^d$ and
$0\le s <\infty$.  With such choices  equation \eqref{eq.v1}
takes the convenient form
\begin{equation}\label{eq.v}
\partial_s v  =\nabla_y\cdot (v^{m-1}\nabla_y v)+
\frac{2}{1-m}\nabla_y\cdot (y v)
        =\nabla_y\cdot\left[v\nabla_y\left(\frac{v^{m-1}-V_D^{m-1}}{m-1}\right)\right]
\end{equation}
This is a  convenient choice since the
stationary states are now given by
\begin{equation}\label{2.7}
V_D(y)= (D+|y|^2)^{-1/(1-m)}, \quad D>0,
\end{equation}
which is just the profile $B_D$ of \eqref{fc} without the undesired
constant $k$. We end this paragraph by noting that for $m=m_*$ the exponent
in the above stationary profile is $-1/(1-m)= -(d-2)/2$, so that $V_D(y)$
decays at infinity like $O(|y|^{-(d-2)})$, i.e., like the
stationary fundamental solution or harmonic potential. This is one
of the reasons that makes $m_*$ special.

%%%%%%%%%%%%%%%%%%%%%%%%%%%%%%%%%%%%%%%%%%%%%%%%%%
\subsection{Stabilization Result}\label{s2.2}

In paper \cite{BBDGV} we have
shown stabilization of solutions of equation \eqref{eq.v} towards
one of the stationary profiles $V_D$ for initial data that are not
very far from $V_D$ to start with. We can write the assumptions on
the {\sl initial conditions\/} in terms of either $u_0$ or $v_0$.
The assumptions on $u_0$ are

\noindent (H1) $u_0$ is a non-negative function in $\LL^1_{\rm
loc}(\RR^d)$ and that there exist positive constants $T$ and
$D_0>D_1$ such that
\begin{equation*}
U_{D_0,T}(0,x)\le u_0(x)\le U_{D_1,T}(0,x)\quad\forall\;x \in\RR^d.
\end{equation*}
\noindent (H2) There exist $D_*\in [D_1,D_0]$ and
$f(|\cdot|)\in\LL^1(\RR^d)$ such that
\begin{equation*}
\big|u_0(x)-U_{D_*,T}(0,x)\big|\le f(|x|)\quad\forall\;x \in\RR^d.
\end{equation*}
In the case $m<m_c$ under consideration here, (H1) implies in particular that the extinction
occurs at time $T$. Moreover, when $m>m_*$ (H2) follows from (H1) since the
difference of two pseudo-Barenblatt solutions is always integrable. For
$m\leq m_*$ this is no more true, and (H2) is an additional restriction.

In terms of $v_0$, conditions (H1) and (H2) can be rewritten as follows.

\noindent (H1') $v_0$ is a non-negative function in $\LL^1_{\rm
loc}(\RR^d)$ and there exist positive constants $D_0> D_1$ such
that
\begin{equation*}
\label{eq:assumptionv}
V_{D_0}(y)\le v_0(y)\le V_{D_1}(y)\quad\forall\;x \in\RR^d.
\end{equation*}
(H2') There exist $D_*\in [D_1,D_0]$ and $\tilde
f(|\cdot|)\in\LL^1(\RR^d)$ such that
\begin{equation*} \label{eq:assumptionmsmallermstarv}
\big|v_0(y)-V_{D_*}(y)\big|\le \tilde f(|y|)\quad\forall\;y \in\RR^d.
\end{equation*}
We point out that condition (H1') means a decay for large $y$ of the form
$$v_0(y)=|y|^{-2/(1-m)}(1-c(y)|y|^{-2})
$$
with $c(y)$ bounded above and below away from zero.
Moreover, (H2') imposes a stronger decay condition for $m\le m_*$. Notice we can take $\tilde f(|y|)=T^{-d\beta}f(|y|/aT^\beta)$, so
that they can be identified up to an elementary scaling.

As a starting point for our asymptotic study, we state the result of \cite{BBDGV}
about the convergence of $v(t)$ towards a unique
Barenblatt profile.

%------------------------------------------------------------------------------
\begin{thm}[Convergence to the asymptotic profile]\label{Thm:A1}
Let $d\ge3$, $m<1$. Consider the solution $v$
of~\eqref{eq.v} with initial data satisfying {\rm
(H1')-(H2')}.
\begin{enumerate}
\item[{\rm (i)}] For any $m>m_*$, there exists a unique
$D_*\in[D_1,D_0]$ such that $\int_{\RR^d}(v(s)-V_{D_*})\dx=0$ for
any $t>0$. Moreover, for any $p\in (q(m),\infty]$,
$\lim_{t\to\infty}\int_{\RR^d}|v(s)-V_{D_*}|^pdy=0$. \item[{\rm
(ii)}] For $m\le m_*$, $v(s)-V_{D_*}$ is integrable,
$\int_{\RR^d}(v(s)-V_{D_*}){\rm d}y=\int_{\RR^d}(v(0)-V_{D_*}){\rm
d}y$ and $v(s)$ converges to $V_{D_*}$ in $\LL^p(\RR^d)$ as
$t\to\infty$, for any $p\in (1,\infty]$. \item[{\rm (iii)}] {\rm
(Convergence in Relative Error)} For any $p\in (d/2,\infty]$,
\begin{equation}\label{CRE}
\lim_{t\to\infty}\left\|{v(s)}/{V_{D_*}}-1\right\|_{p}=\;0\;.
\end{equation}
\end{enumerate}
\end{thm}
%------------------------------------------------------------------------------
For simplicity, we  write $v(s)$ instead
of $y \mapsto v(s,y)$ whenever we want to emphasize the dependence
on the time $s$. The exponent $q(m)$ is defined as the infimum of all positive real numbers $p$ for which
two Barenblatt profiles $V_{D_1}$ and $V_{D_2}$ are such that
$|V_{D_1}-V_{D_2}|$ belongs~to~$\LL^p(\RR^d)$:
\begin{equation*}\label{q_0}
q(m):=\dfrac{d(1-m)}{2(2-m)}\;.
\end{equation*}
We see that $q(m)>1$ if $m\in(0,m_*)$, $q(m_*)=1$, and $q(m)<1$ if
$m>m_*$. In case $m>m_*$, the value of $D_*$ can be computed at
$s=0$ as a consequence of the mass balance law
$\int_{\RR^d}(v_0-V_{D_*})\dx=0$, and then the conservation result
holds for all $s>0$ as is proved in the paper \cite{BBDGV}. On the
other hand, in the case $m\le m_*$ the mass balance does not make
sense, but $D_*$ is determined by Assumption (H2'). In this case, the
presence of a perturbation of $V_{D_*}$ with nonzero mass, does not
affect the asymptotic behavior of the solution to first order.

%%%%%%%%%%%%%%%%%%%%%%%%%%%%%%%%%%%%%%%%%%%%%%%%%%%%%%%%%%%%%

\subsection{Relative error, entropies, and linearization}
\label{entropies}

The deeper stabilization analysis of equation \eqref{eq.v}
leads to an interesting connection with a family of Poincar{\'e}-Hardy
functional inequalities. In  this way, we obtain stabilization
rates that are exponential in the new time $s $, which means that
they are power-like in the original time. The exponent $m_*$
appears precisely as the only exponent for which the linearized
analysis based on Poincar{\'e}-Hardy inequalities fails and the
corresponding rates are not obtained by that method. We shall prove below that the
linearized analysis when $m$ takes the special value $m_*$ leads to a different functional framework and
the actual rates are different, and actually slower.

In any case, the approach and the use of entropies
starts in the same way. Let $v$ be a solution to the rescaled Fokker-Plank equation
\eqref{eq.v}, and let $V_*=V_{D_*}$ be the Barenblatt solution
mentioned in Theorem \ref{Thm:A1}. We pass to the quotient
$w(s,y)=v(s,y)/V_*(y)$. Notice that $w-1=(v-V_*)/V_*$ is the
relative error of $v$ with respect to $V_*$.  Notice also that, by
straightforward calculations, our running assumptions imply that
$W_0\le w\le W_1$, where \ $W_0=\left(D_*/D_0\right)^{1/(1-m)}<1$
and $W_1=\left(D_*/D_1\right)^{1/(1-m)}>1$.

The equation for $w$ reads
\begin{equation}\label{eq.w}
\partial_s w=\frac{1}{V_*}\nabla\cdot\left[wV_*\nabla\left(\frac{w^{m-1}-1}{m-1}V_*^{m-1}\right)\right]
\end{equation}
In terms of $w$, we define the {\sl relative entropy\/}
\begin{equation}\label{Entropy.Quotients}
{\mathcal F}[w]:= \frac1{1-m}\int_{\RR^d}\left[(w-1)-\frac
1m(w^m-1)\right]V_{*}^m\dy
\end{equation}
Strictly speaking, we are assuming that a time $s\ge 0$ is given
and then we get ${\mathcal F}(w(s))$. In terms of $v$, when $m$ is
sufficiently close to one,  it can be derived as $E(v)-E(V_*)$
where
\begin{equation}E[v]=:\frac1{1-m}\int_{\RR^d}\left[ vV_*^{m-1}-\frac1m
v^m\right]\dy
 \end{equation}
(for $m$ farther away from 1 both $E[v]$ and $E[V_*]$
become infinite and only the expression for the difference, the
relative entropy, makes sense). We also introduce
the {\sl relative Fisher information\/}
\begin{equation}\label{Fisher.Quotients}
\mathcal{I}[w]
    =\int_{\RR^d}\left|\nabla\left(\frac{w^{m-1}-1}{m-1}V_*^{m-1}\right)\right|^2 V_* w \dy
    =\int_{\RR^d}\left|\nabla\left(\frac{v^{m-1}-V_*^{m-1}}{m-1}\right)\right|^2
v \dy
\end{equation}
(again, we should have written $\mathcal{I}[w(s)]$). By differentiation in time and using the equation, we get
\begin{equation}\label{entropydiff}
\frac{{\rm d}{\mathcal F}[w(s)]}{{\rm d}s}=-{\cal I}[w(s)]\,, \qquad\forall s>0\,.
\end{equation}
For a detailed proof of this time derivation, we refer to Proposition 2.6 of \cite{BBDGV}.

\medskip

We now introduce the linearization idea in \cite{BBDGV} that allows to treat the long-time behaviour
of $w$. It consists in writing the relative error in the form
\begin{equation}\label{2.14}
w(s,y)-1=\varepsilon g(s,y)V_*^{1-m}(y)
\end{equation}
where the choice of weight $V_*^{1-m}$ is crucial. After a brief formal computation
we obtain the differential equation for $g$ that is implied by \eqref{eq.w} in the limit
$\varepsilon\to 0$:
\begin{equation} \label{lin.g}
\partial_s g =V_*^{m-2}\nabla\cdot\left(V_*\nabla g\right).
\end{equation}
Actually, since Theorem \ref{Thm:A1}, formula \eqref{CRE}, implies
that $w\to 1$ as $s\to\infty$, the factor $\varepsilon$  will not
be needed in the actual linearization step.

Our next task is to study this linear flow; then, we shall have to
relate the actual nonlinear flow to its linearized approximation.
  But let us point out that we will not need to prove
the convergence of solutions of the original problem to solutions
of the linear problem, the analysis is rather based on the
relationship between the two linear quantities, entropy and Fisher
information, associated to equation \eqref{lin.g}, and the close
similarity of these linear  quantities and the previously defined
nonlinear ones. These facts plus \eqref{entropydiff} produce the
desired convergence result.

In the cases $m<1$, $m\ne m_*$, a suitable functional setting was
found where the functional inequalities of Hardy-Poincar\'e type
corresponding to the linear flow implied the existence of a
spectral gap. According to more or less standard theory, existence
of such a gap implies exponential decay rates (in $s$) of the
norms and entropy of the solutions of the linear flow. A delicate
analysis of comparison of entropy and Fisher information between
the linear and nonlinear flow allowed finally to transfer the
result about decay rate to the original nonlinear flow. See full
details in \cite{BBDGV}.

The problem arising when $m=m_*$ is the absence of spectral gap. We
shall prove below that this is essential, in fact the actual rates
are not exponential but power--like in $s$. This is related to the
heat kernel behaviour of the operator appearing in \eqref{lin.g},
$V_*^{m-2}\nabla\cdot\left(V_*\nabla g\right)$, acting in a
suitable weighted Hilbert space.
Details will be given in Section \ref{sect-lin}, where a sharp
power-like decay for the heat kernel is proved using a most fortunate coincidence, i.\,e., the
representation of the linear semigroup  as the heat flow on a
conformally flat Riemannian manifold. It will also be proved that
no Hardy--type inequality can hold for the quadratic form
associated to the generator, so that it is hopeless to use the
same line of reasoning of paper \cite{BBDGV}.

%%%%%%%%%%%%%%%%%%%%%%%%%%%%%%%%%%%%%%%%%%%%%%%%%%%%%%%%%%%%%%%%%%%%%%%%%%
\section{Statement of the main results for $m=m_*$}
\label{ssec.statement}

We are now ready to state our main results. We use the notations
$v(s), w(s)$ instead of $v(s,y), w(s,y)$ and $u(t)$ instead of
$u(t,x)$ when the dependence on time is stressed.

We prove convergence of $v(s)$ to the appropriate Barenblatt
profile in several senses. More precisely we prove quantitative
bounds on the convergence in suitable L$^p$ norms, on the
convergence of moments, and on the uniform convergence of all
derivatives. Convergence takes place {\it with the same rate of
the linearized case}.

\begin{thm}[Convergence with rate to the asymptotic profile]
\label{thm.conv.entropy} Consider a solution $v$ of the equation \eqref{eq.v}
such that $v_0$ satisfies {\rm (H1')-(H2')} and fix some $s_0>0$.
Then, the entropy of the quotient variable satisfies
\begin{equation}
\mathcal{F}[w(s)]\leq K s^{-1/2}\quad\forall\;s\geq s_0\;.
\end{equation}
 for some $K=K(v_0,s_0)$. As a consequence, for any $\vartheta\in\left[0,\frac d2\right]$, there
exists a positive constant $K_\vartheta$ such that
\begin{equation}
\left\||y|^\vartheta(v(s)-V_{D_*})\right\|_{2}\leq K_\vartheta
s^{-1/4}\quad\forall\;s\geq s_0\;.
\end{equation}
\end{thm}

The analysis of the linearized equation indicates
that this rate should be optimal. We also have convergence without
weights in suitable $\LL^p$ and $C^j$ spaces with the same rates,
where we use interior regularity theory for parabolic equations:

\begin{cor} \label{Conv.Weight}\ {\rm (i)} For any $q\in(1,\infty]$, there exists a positive
constant $K(q)$ such that
\begin{equation}
\|v(s)-V_{D_*}\|_q\le K(q)s^{-1/4}\quad\forall\;s\geq s_0\;.
\end{equation}
\noindent {\rm (ii)} For any $j\in\mathbb{N}$ there exists a
positive constant $H_j$ such that
\begin{equation}
\|v(s)-V_{D_*}\|_{C^j(\RR^d)}\le H_js^{-1/4}\quad\forall\;s\geq s_0\;.
\end{equation}
\end{cor}

These power-decay results are in contrast with the exponential rates
obtained in \cite{BBDGV} for $-\infty < m <1$ and $m \ne m_*$.
Rescaling back to the original space--time variables one gets the
following result  which can be called {\sl intermediate
asymptotics}.

\begin{cor}\label{Cor:A2}
Consider a solution $u$ of~\eqref{1.1} with  $m=m_*$, with initial data satisfying {\rm
(H1)-(H2)}, and extinction time $T$. For $t$ sufficiently close to $T$ and for any $q \in (1,\infty]$,
there exists a positive constant $C$ such that:
\begin{equation*}
\|u(t)-U_{D_*}(t)\|_q\le C\,(T-t)^{\sigma(q)}\,\log\left(T/(T-t)\right)^{-1/4}.
\end{equation*}
with $\sigma(\infty)=d(d-2)/4$, and $\sigma(q)=\sigma(\infty)(q-1)/q$ \ for $q<\infty$. \end{cor}

\noindent We also obtain a quantitative bound on the decay of the
{\it relative error of $v(s)$ with respect to  $V_{D_*}$}.

\begin{cor}[Decay of Relative Error]\label{thm:CRE-exp}
Consider a solution $v$ of \eqref{eq.v} such that $v_0$ satisfies {\rm (H1')-(H2')} and fix some $s_0>0$. Then for any $q\in(d/2,\infty]$ and all $\varepsilon>0$ there exists a positive constant
$\mathcal{C}_q$ such that
\begin{equation}
\big\|{v(s)}/{V_{D_*}}-1\big\|_q\le \mathcal{C}_qs^{-\frac {1-\varepsilon}{d}}\quad\forall\;s\geq s_0\;.
\end{equation}
If $q=d/2$ there is a positive constant $\mathcal{C}$ such that
\begin{equation}
\big\|{v(s)}/{V_{D_*}}-1\big\|_{d/2}\le \mathcal{C}s^{-\frac 1{d}}\quad\forall\;s\geq s_0\;.
\end{equation}
Finally we also have, for all $j\in{\mathbb N}$, that there exists
a positive constant $\mathcal{C}_j$
\begin{equation}
\big\|{v(s)}/{V_{D_*}}-1\big\|_{C^j(\RR^d)}\le \mathcal{C}s^{-\frac {1-\varepsilon}{d}}\quad\forall\;s\geq s_0\;.
\end{equation}
\end{cor}
Notice that, besides having a quantitative bound, we have some other
improvements on Theorem \ref{Thm:A1} first because the value $q=d/2$
is now allowed and because convergence of $C^j$ norms is also dealt
with. The constants involved depend also on $m,d,D_0, D_*, D_1$, but
also on the solution at time $s_0$ through the relative mass
(conserved along the evolution) and through the uniform bound $c_0$
on the ratio $\int_{{\mathbb R}^d} |\nabla
v(y)|^2V_D(y)\dy/\|v\|_1^2\le c_0$\,.

%%%%%%%%%%%%%%%%%%%%%%%%%%%%%%%%%%%%%%%%%%%%%%%%%%%%%%%%

\section{Analysis of the linear case}
\label{sect-lin}

We address now a central topic of the paper, i.e., establishing of the
long-time behaviour of the linearized flow in the still open case with exponent
$m_*=(d-4)/(d-2)$. The clue to our study of the linearized flow in
this case is to interpret it as the heat flow of the
Laplace-Beltrami operator of a suitable Riemannian  manifold
$(M,{\bf g})$, with a metric ${\bf g}$ which is conformal to the standard
$\RR^d$ metric. Studying the pointwise heat kernel behaviour allows
to prove Nash and log-Sobolev inequalities associated to the generator.
Such inequalities will later on allow us to study the nonlinear
evolution as well, and to determine its asymptotics, which will be shown to proceed
with the same rate of convergence as the linearized one. Since the
study can have independent interest,  we replace $g$ by $v$,  $y$ by
$x$, and $s$ by $t$  throughout the section to conform to more
standard notations.

\subsection{Linear equation and geometry}\label{linear}
Given $m<1$ and $D>0$, we consider the operator given on
$C_c^\infty({\mathbb R}^d)$ ($d\ge3$) by
\beq\label{op}
L_mv=
(D+|x|^2)^{(2-m)/(1-m)}\nabla\cdot\left(\frac{\nabla
v}{(D+|x|^2)^{1/(1-m)}}\right)= V_D^{m-2}\nabla\cdot\left(V_D\nabla
v\right). \eeq
We recall that for $m=m_*$ the following holds:  $1/(1-m)=(d-2)/2$, and $V_D^{m-2}(x)=(1+|x|^2)^{d/2}$.
We have dropped the index $*$ from $D_*$
to simplify the notation, since the particular value of $D$ has no role here.
 We shall think of this operator as acting on the
Hilbert space $H_m={\rm L}^2({\mathbb R}^d, \rd\mu)$ with
$d\mu=V_D^{2-m}\rd x$. To define it more precisely we construct the
quadratic form
\beq\label{form}
I_m[v]=\int_{{\mathbb R}^d}
\frac{|\nabla v(x)|^2}{(D+|x|^2)^{1/(1-m)}}{\rm d}x=\int_{{\mathbb
R}^d} |\nabla v(x)|^2V_D(x){\rm d}x,\ \ \ u\in C^\infty_c({\mathbb
R}^d).
\eeq
Then, ${I}_m$ is closable in $H_m$ (for a quite general
result implying the validity of the above assertion see e.g.
\cite{D}, Section~4.7).  We denote again by $-L_m$ the unique
nonnegative self--adjoint operator in $H_m$ associated with its
closure. In fact $L_m$ has the above explicit expression \eqref{op} on smooth
compactly supported functions. There is a particular value of $m$
for which the above operator can be seen as the Laplace-Beltrami
operator of a certain Riemannian manifold $(M,{\bf g})$, as we shall
show. This in particular will imply (since $M$ turns out to be
complete) that $L_m$ is essentially self--adjoint on
$C_c^{\infty}(M)$ by a result of Calabi (see e.g. \cite{D}, Theorem
5.2.3). Consider indeed the following manifold, denoted by $M$,
given by ${\mathbb R}^d$ endowed with the Riemannian, conformally
flat metric defined, in Euclidean (global) coordinates, by
\begin{equation}\label{metric}
{\bf g}_\alpha(x)=(D+|x|^2)^{-\alpha}{\bf I},
\end{equation}
where ${\bf I}$ is the Euclidean metric and $|\cdot|$ is the
Euclidean norm. We denote by $\mu_{{\bf g}_\alpha}$ the Riemannian
measure, by $|{\bf g}_\alpha|=\mbox{\rm det}({\bf g}_\alpha)$
the determinant of the metric tensor,
by $\nabla_\alpha$ the Riemannian gradient and by $\Delta_\alpha$
the Laplace-Beltrami operator, defined on L$^2(\mu_{{\bf
g}_\alpha})$, associated to the given metric.

\begin{lem}
The Laplace-Beltrami operator $\Delta_{\alpha}$
coincides with $L_m$, precisely when $\alpha=1$ and $m=m_*:=(d-4)/(d-2)$,
both as concerns its explicit expression (in Euclidean
coordinates) and as concerns the Hilbert space it acts on.
\end{lem}
\proof We notice that for the above choice of metric we have
 \[ \sqrt{|{\bf g}_\alpha|(x)}=(D+|x|^2)^{-\alpha
d/2},\ \ \ g_\alpha^{ij}(x)=(D+|x|^2)^{\alpha}\delta^{ij}.
\]
Then we have that the Dirichlet form associated to $\Delta_{\alpha}$ is
given, on test functions, by \beq
\begin{aligned}J_\alpha(v)&:=\int_M {\bf g}_\alpha(\nabla_\alpha v,\nabla_\alpha
v)\rd\mu_{{\bf g}_\alpha}=\int_{{\mathbb R}^d} \sqrt{|{\bf
g}_\alpha|(x)}g^{ij}_\alpha(x)\frac{\partial v}{\partial
x^i}\frac{\partial
v}{\partial x^j}\dx\\
&=\int_{{\mathbb
R}^d}(D+|x|^2)^{(-d\alpha/2)+\alpha}|\nabla_ev(x)|^2\dx
\end{aligned}
\eeq where $\nabla_e$ is the Euclidean gradient and the summation
convention is used. Then we notice that the conditions that identify $\Delta_\alpha$ with
$L_m$:
\[\begin{aligned}
&\sqrt{|{\bf g}_\alpha|(x)}=(D+|x|^2)^{-(2-m)/(1-m)}\\
&\sqrt{|{\bf
g}_\alpha|(x)}g^{ij}_\alpha(x)=(D+|x|^2)^{-1/(1-m)}\delta^{ij}
\end{aligned}
\]
force $\alpha, m$ to be related by $(d\alpha/2)-\alpha=1/(1-m)$ and
$d\alpha/2=(2-m)/(1-m)$. This is equivalent to $\alpha=1$,
$m=(d-4)/(d-2)=m_*$ as claimed.\qed

\medskip

We shall now compute, in the case discussed in the above Lemma, the
Ricci curvature of $(M,{\bf g}_\alpha)$. Hereafter we shall drop the index $\alpha$, since we always choose $\alpha=1$. We put $D=1$ for simplicity without loss of generality.

\begin{lem}\label{lemma.Rij}
Then the Ricci curvature of $(M,{\bf
g}_{\alpha=1})$ is given, in Euclidean coordinates, by
\beq
{R}_{ij}=-\frac{(d-2)x_ix_j}{(1+|x|^2)^{2}}+
\left[\frac{(d-2)|x|^2+ 2(d-1) }{(1+|x|^2)^{2}}\right]\delta_{ij},
\eeq
where we write {\rm Ric}\,$=(R_{ij})$. In particular {\rm
Ric}\,$>0$ on $M$, such lower bound cannot be improved, and {\rm
Ric} is bounded on $M$. Actually, $R_{ij}(x)=O(|x|^{-2})$ as $|x|\to\infty$ in the transversal directions and it behaves as  $O(|x|^{-4})$ in the radial directions. Finally, the scalar curvature is given by
\begin{equation}
R=(d-1)\frac{2d+ (d-2)|x|^2}{1+|x|^2}.
\end{equation}
\end{lem}

We will postpone the proof of these formulas to appendix A1 not to break the flow
of the exposition. It immediately follows that the  symmetric tensor Ric is positive;
indeed, given $\xi\in {\mathbb R}^d$, we have
\[
R_{ij}(x)\xi_i\xi_j
\ge\frac{2(d-1)}{(1+|x|^2)^2}|\xi|^2>0.
\]
 The boundedness of Ric is clear from its explicit
expression. Note that for $d=2$ we are dealing with an Einstein metric, $\textrm{Ric}=k\, {\bf g}$ (actually,
it is Hamilton's cigar soliton to the Ricci flow, \cite{Ham, ChowK}),
but for $d\ge 3$ it is not.

\medskip

Let us continue with the asymptotic analysis of the flow. By a celebrated result of Li and Yau \cite{LY}, the heat kernel $K(s,x,y)$ of the Laplace--Beltrami operator of a complete
Riemannian manifold $(M,{\bf g})$ with nonnegative Ricci curvature is
pointwise comparable with the quantity
\[
\frac1{{\rm Vol}[B(x,\sqrt t)]}e^{-c\frac{d^2(x,y)}{t}}
\]
where $d(\cdot\,,\cdot)$ is the Riemannian distance in $(M,{\bf g})$, $B(x,r)$ is the
Riemannian ball centered at $x$ and of radius $r$ and Vol is the
Riemannian volume. More precisely,
\begin{cor}
For all small positive $\varepsilon$ there exists positive constants $c_1,c_2$ such that
\[
\frac{c_1(\varepsilon)}{{\rm Vol}[B(x,\sqrt t)]}e^{-\frac{d^2(x,y)}{(4-\varepsilon)t}}\le
K(t,x,y)\le \frac{c_2(\varepsilon)}{{\rm Vol}[B(x,\sqrt
t)]}e^{-\frac{d^2(x,y)}{(4+\varepsilon)t}}
\]
for all $x,y\in M$, $t>0$.
\end{cor}
We recall that the Li--Yau bounds require completeness, a property
which clearly holds for the manifold we are considering. We use
the notation $a\wedge b=\min\{a,b\}$.

\begin{cor}
The heat kernel satisfies the following properties:
\begin{equation}
\label{on-diag}\begin{split}
&K(t,x,x)\appros_{t\to0}\left(1\wedge\frac1{|x|}\right)\frac1{t^{\frac d2}},\\
&K(t,x,x)\le%\left(1\wedge\frac1{|x|}\right)
\frac C{t^{\frac 12}}\ \ \forall t\ge1, \forall x\in{\mathbb R}^d,
\end{split}
\end{equation}
where $f_1\appros\limits_{t\to t_0} f_2$ means that there exists two constants $c_1,c_2>0$ such that $c_1f_1\le f\le c_2f_2$ near $t_0$.
\end{cor}
\proof
First notice that
\[d(0,x)=\int_0^{|x|}\frac1{\sqrt{1+t^2}}\rd t\]
where $|x|$ is the Euclidean length, so that $d(0,x)\sim\log |x|$ for large $|x|$. Hence,
\[\begin{aligned}
\vol(B(0,R))&=\int_{B(0,R)}\sqrt{|{\bf g}|}\rd x=\int_{d(0,x)<R}\frac1{(D+|x|^2)^{d/2}}\rd x\\
&\asint_{R\to+\infty}c\int_{r<e^R}\frac{r^{d-1}}{(D+r^2)^{d/2}}\rd r\asint_{R\to+\infty}cR.
\end{aligned}
\]
Proceeding similarly, one shows that $d(x_0,x)\appros \log\frac{|x|}{|x_0|}$ for large $|x|$ and hence that $\vol(B(x_0,R))\sim c(R+\log|x_0|)$ for large $R$ and, say, $|x_0|\ge2$. The short time behaviour is clearly locally Euclidean, with a weight depending on $x$ given by definition by $1/\sqrt{D+|x|^2}$.
\qed
\noindent{\bf Remark}. The above corollary extends, for the present choice of the parameter $m$, the result of \cite{D}, Th. 4.7.5, in several respects. In fact, in the quoted Theorem the bounds on the heat kernel are from above and for short time only. Notice that the short time bound in the following results matches with the one of \cite{D}. One may notice that, in fact, we have proved the bound $$K(t,x,x)\appros
\frac 1{t^{\frac 12}+\log(1+|x|)}\ \ \ \forall t\ge1, \forall x\in{\mathbb R}^d,$$ although we shall make no further use of it.

\begin{cor}
Each solution to the linear evolution equation $\partial_t v=L_{m_*}v$ corresponding to an initial datum in ${\rm L}^1({\mathbb R}^d,(D+|x|^2)^{(m-2)/(1-m)})$ satisfies the bound
\begin{equation}\label{Ultra.1}
\|v(t)\|_\infty\le H(t)\|v_0\|_1=
\left\{\begin{array}{lll}
c_1\dfrac{\|v_0\|_1}{t^{d/2}} &\mbox{for any }0<t\le 1\\[3mm]
c_2\dfrac{\|v_0\|_1}{t^{1/2}} &\mbox{for any }t>1
\end{array}\right.
\end{equation}
where $c_i$ are positive constants. The power of $t$ cannot be improved for such general initial data, as can be seen by considering the time evolution of a Dirac delta.
\end{cor}

\noindent  {\bf Warning:} Here, the symbol $\|\cdot\|_p$ denotes the
norm in ${\rm L}^p({\mathbb R}^d,\rd \mu_*)$, where
$\rd\mu_*=V_{D}^{2-m_*}(x)\rd x$, and we know that $V_{D}^{2-m_*}(x)=(D+|x|^2)^{-d/2}$.  This notation will be kept in
the next three sections.

%%%%%%%%%%%%%%%%%%%%%%%%%%%%%%%%%%%%%%%%%%%
\subsection{Functional Inequalities}\label{funct}

We recall that $I_{m_*}[v]=\int_{{\mathbb R}^d} |\nabla v(x)|^2
V_{D}\,{\rm d}x$ on  smooth compactly supported functions. The
domain of its closure will be indicated by ${\rm Dom}\,(I_{m_*})$.

\begin{cor}
There is a family of logarithmic Sobolev inequalities
\begin{equation}
\int_{{\mathbb R}^d} v^2\log\left(\frac{v}{\|v\|_2}\right)\rd \mu_*
\le \varepsilon I_{m_*}[v]+\beta(\varepsilon)\|v\|_2^2
\end{equation}
valid for all $v\in {\rm Dom}(I_{m_*})\cap {\rm L}^1({\mathbb
R}^d,\rd \mu_*)\cap {\rm L}^\infty({\mathbb R}^d,\rd \mu_*)$ and all
positive $\varepsilon$, where $\beta(\varepsilon)=c-\frac
d4\log\varepsilon$ for $\varepsilon<1$, $\beta(\varepsilon)=c-\frac
14\log\varepsilon$ for $\varepsilon\ge1$, and $c$ is a suitable
positive constant.
\end{cor}

\proof We have $\|v(s)\|_\infty\le Cs^{-1/2}\|v_0\|_1$ for large
$s$. Interpolating between such bound and the L$^\infty$
contractivity property (valid since $I_{m_*}$ is a Dirichlet form)
shows that $\|v(s)\|_\infty\le Cs^{-1/4}\|v_0\|_2$ for large $s$.
Similarly, $\|v(s)\|_\infty\le Cs^{-d/4}\|v_0\|_2$ for small $s$.
The validity of such ultra-contractive bounds for the solution of
the linear evolution considered is known to be equivalent, by
\cite{D}, Example 2.3.2, to the stated logarithmic Sobolev
inequalities for the initial datum $u_0$ if it belongs to ${\rm
Dom}(I_{m_*})\cap {\rm L}^1({\mathbb R}^d,\rd \mu_*)\cap {\rm
L}^\infty({\mathbb R}^d,\rd \mu_*)$. At this point the evolution has
no role anymore and to avoid confusions we choose to write $v$
instead of $u_0$ in the statement. \qed

The next consequences we draw involve the recurrence of the
semigroup considered.
\begin{cor}
The semigroup $\{T_s\}_{s\ge0}$ associated to $L_{m_*}$ is
recurrent. In particular, $L_{m_*}$ does not admit a (minimal) positive Green
function and the manifold $({\mathbb R}^d,{\bf g}_{\alpha=-1})$ is
parabolic.
\end{cor}

\proof It suffices to note that a semigroup $\{T_s\}_{s\ge0}$ is, by
definition, transient, iff $\int_0^\infty T_sv\rd s$ is a.e. finite
for all $v\in\LL^2({\mathbb R}^d,\rd \mu_*)$. This of course does
not hold in the present case because of the $s^{-1/2}$ behaviour for
long times of the heat kernel. \qed
\begin{cor}
There is no bounded, strictly positive, $\mu_*$--integrable function
$h$ such that
\[
\int_{{\mathbb R}^d}|v|h\rd \mu_* \le I_{m_*}[v]^{1/2}
\]
for all $v\in {\rm Dom}\,(I_{m_*})$.
\end{cor}
\proof The existence of a function $h$ with the stated properties
is equivalent to the transience of the semigroup at hand, by
\cite{F}, Th. 1.5.1. \qed

\begin{cor}\label{No.Hardy.No.Party}
There is no bounded, strictly positive, $\mu_*$-integrable function
$h$ such that for all $v\in {\rm Dom}\,(I_{m_*})$
\[ \int_{{\mathbb
R}^d}v^2h\rd \mu_* \le I_{m_*}[v].
\]
\end{cor}

\proof Since $h$ is assumed to be integrable so that $h\rd \mu_*$ is
a finite measure, that we can normalize to 1, one would have by
H\"older inequality that
\[
\int_{{\mathbb R}^d}|v|h\rd \mu_* \le \left(\int_{{\mathbb
R}^d}v^2h\rd \mu_*\right)^{1/2}\le (I_{m_*}[v])^{1/2}
\]
for all $v\in {\rm Dom}(I_{m_*})$, contradicting the above result.
\qed

\noindent{\bf Remark}. The above results prove that Hardy--type
inequalities relative to the Dirichlet form $I_{m_*}$ and to a
strictly positive integrable weight $h$ cannot hold, even if $h$ is
required to be bounded. This shows that the strategy of
\cite{BBDGV}, which relied heavily on the validity of Hardy--type
inequalities and allowed to deal with the case $m\not=m_*$ cannot be
adapted to the present situation.

The ultra-contractive bounds discussed above can also be related to
the validity of Nash inequalities for $I_{m_*}$. In fact we prove
now some inequalities of that type in weighted Sobolev spaces which
will be very important when dealing with the nonlinear evolution.
Such inequalities play here the role that Hardy--type inequalities
played in the case $m\neq m_*$ studied in \cite{BBDGV}, cf. also Section \ref{m.neq.mstar}. The
following crucial result is a purely functional inequality which is
proved using the linear evolution only, but will turn out to be the
key point for the study of the nonlinear evolution as well.

\begin{prop}\label{GNprop}
For all $v$ such that $I_{m_*}[v]/\|v\|_1^2\le c_0$ for some
$c_0>0$, the following Gagliardo--Nirenberg inequality
holds true:
\begin{equation}\label{GN}
\|v\|_2^2\le KI_{m_*}[v]^{1/3}\|v\|_1^{4/3},
\end{equation}
for all $v\in \LL^2({\RR}^d,\rd \mu_*)\cap{\rm Dom}(I_{m_*})$,
where the positive constant $K$ depends on $c_0$, and
diverges as $c_0\to+\infty$.
\end{prop}
\proof\; To get the claim, first interpolate between the bound
$\|v(s)\|_\infty\le H(s)\|v_0\|_1$ and the L$^1$ contraction
property to get $\|v(s)\|_2\le H(s)^{1/2}\|v_0\|_1$. From this
starting point we can use a known argument, cf.  \cite{D}, and we
briefly recall it for the sake of completeness. In fact, use the
semigroup property and the fact that $I_{m_*}[v(s)]$ is nonincreasing as
a function of $s$ to write
\[
\begin{aligned}
H(s)\|v_0\|_1^2&\ge (v(s), v(s))=(v(2s),v_0)\\
&=(v_0,v_0)-\int_0^{2s}I_{m_*}\left[v\left(\frac \lambda2\right)\right]\rd \lambda\\
&\ge (v_0,v_0)-2sI_{m_*}[v_0].
\end{aligned}
\]
Therefore,
\begin{equation}\label{optimize}
\|v_0\|_2^2\le 2sI_{m_*}[v_0]+H(s)\|v_0\|_1^2.
\end{equation}
It would then be easy to minimize the r.h.s. of the latter formula
should one have $H(s)=cs^{-\alpha}$ for all $s>0$. The fact that
$H(s)$ has such form with different powers of time when $s$ is
small and when $s$ is large forces us to proceed as follows.
Assuming that $I_{m_*}[v_0]$ and $\|v_0\|_1$ are not zero, the
right hand side takes the value infinity both as $s\to 0$ and
$s\to\infty$ hence there is a minimum for one or several
intermediate values of $s$. We want to take a particular value of
$s$ that almost minimizes the above formula, and we want that
value to correspond to the range of not small $s$ where
$H(s)=c_2s^{-1/2}$. Since we assumed that
$I_{m_*}[v_0]/\|v_0\|_1^2\le c_0$ for some $c_0>0$, we consider
the 1-parameter quantity
\[
s_\alpha=\alpha\left[\frac{\|v_0\|_1^2}{I_{m_*}[v_0]}\right]^{2/3}
\]
and observe that, trivially,
\[
s_\alpha> 1\qquad\iff\qquad \alpha >c_0^{2/3}.
\]
We choose $\alpha$ accordingly (so that it is bounded away from
zero) and plug the corresponding $s_\alpha$ into \eqref{optimize},
noticing that for $s>1$ we have $H(s)=c_2s^{-1/2}$; with these choices \eqref{optimize} becomes
\begin{equation}
\|v_0\|_2^2\le K\|v_0\|_1^{4/3}I_{m_*}[v_0]^{1/3}
\end{equation}
with $K=K(\alpha)=\alpha+c_2\alpha^{-1/2}$. This concludes the
proof.\qed

%%%%%%%%%%%%%%%%%%%%%%%%%%%%%%%%%%%%%%%%%%%%%%%%%%%%%%%%%%%%%%
\subsection{Mass conservation for the linear flow}
We introduce here the calculation of ``conservation of mass'' for
the linear semigroup. As usual we put $V^{2-m_*}_{D}\rd x=\rd
\mu_*$.

\begin{lem} The following property of mass conservation holds true
for every {\it nonnegative} $v\in \LL^1(\rd \mu_*)$:
\begin{equation}
\frac{\rd}{\dt }\int v \rd \mu_*=0.
\end{equation}
\end{lem}

\noindent We give two proofs of the result, first a quite direct
one and then a proof relying on the special geometric nature of
the linear flow.

\noindent{\bf First proof}. We use the specific form of the weights
involved for a direct calculation, first for an initial
datum belonging to ${\rm L}^1(\rm{d}\mu_*)\cap {\rm
L}^2(\rm{d}\mu_*)$. Choose a test function $\varphi_R$
supported in the Euclidean ball $B_{2R}$ with
$\varphi_R=1$ on $B_R(0)$. Let $t\ge t_0>0$ and compute, for any such $t$:
\begin{equation}\begin{split}
&\left|\frac{\rd}{\dt }\int_{\RR^d} v\varphi_R\rd \mu_* \right|
    =\left|-\int_{R\le |x|\le 2R} \nabla v\cdot \nabla\varphi_R V_{D}\rd x\right|\\
    &\le \int_{R\le |x|\le 2R} |\nabla V_D|\,|\nabla\varphi_R| v \rd x
     + \int_{R\le |x|\le 2R} |\Delta \varphi_R|\,v\,V_D\rd x
    \le \frac{c(m,d)}{R^2}\int_{R\le |x|\le 2R} v\,V_D\rd x
\end{split}
\end{equation}
because $|\nabla V_D|\le c_0(m,d) V_D/R$ since it is not restrictive to assume that $|\nabla \varphi_R|\le c_1/R$ and $|\Delta \varphi_R|\le c_2 /R^2$ whenever $R\le |x|\le 2R$. By H\"older inequality we obtain that
\[
\int_{R\le |x|\le 2R} v\,V_D\rd x
    \le \left(\int_{R\le |x|\le 2R} v^2\,V_D^{2-m_*}\rd x\right)^{1/2}
        \left(\int_{R\le |x|\le 2R} V_D^{m_*}\rd x\right)^{1/2}
    \le \varepsilon_R\, R^2
\]
since we let $\varepsilon_R:=\left(\int_{R\le |x|\le 2R} v^2\,V_D^{2-m_*}\rd x\right)^{1/2}$ and it is easy to check that $\int_{R\le |x|\le 2R} V_D^{m_*}\rd x\le c_1 R^4$.
We obtained that
\[
\left|\frac{\rd}{\dt }\int_{\RR^d} v\varphi_R\rd \mu_* \right|\le c_1\varepsilon_R
\]
and we notice that $\varepsilon_R\to 0$ as $R\to \infty$, a fact which holds because $v\in\LL^2(\rd\mu_*)$. This proves that
\begin{equation}
   \lim_{R\to\infty}\left|\int_{\RR^d} v(t_1)\varphi_R\rd \mu_*-\int_{\RR^d} v(t_0)\varphi_R\rd \mu_*\right|\le c_1\lim_{R\to\infty}\varepsilon_R (t_1-t_0)=0
\end{equation}
for any $0\le t_0\le t_1$. We can use dominated convergence in the left-hand side, since the Markov property implies $v(t)\in\LL^1(\rd\mu_*)$ for all $t\ge 0$. This yields the claim for strictly positive times and for initial data belonging to ${\rm L}^1(\rm{d}\mu_*)\cap {\rm
L}^2(\rm{d}\mu_*)$. We can then reach $t=0$ using the strong continuity in $\LL^1(\rd\mu_*)$ of the evolution, and consider general data in $\LL^1(\rd\mu_*)$ by approximation.

\noindent{\bf Second proof}. We can also use a general
argument involving conservation of probability on manifolds with
curvature bounded below. Let $\{T_t\}_{t\ge0}$ be the semigroup
associated to the Laplace--Beltrami operator of the manifold
considered. Then $\{T_t\}_{t\ge0}$ is a Markov semigroup and in
particular it acts on all L$^p$ spaces ($p\in[1,+\infty]$), it is
contractive on any such space and it preserves positivity. We have
shown that the Ricci curvature of $M$ is bounded. An application of
\cite{D}, Theorem 5.2.6 then shows that $\{T_t\}_{t\ge0}$ preserves
the identity: $T_t1=1$. From this, conservation of the L$^1$ norm
for data $v\ge0$ follows. In fact, with the notation $v(t)=T_tv$ and
using the fact that the adjoint of $T_t$ when seen as acting on
L$^1$ is $T_t$ itself but seen as acting on L$^\infty$, we have:
\[\begin{aligned}\displaystyle
\|v(t)\|_1&=\sup_{h\in{\rm
L}^\infty\hskip-3pt,|h|\le1}\left|\int_{{\mathbb R}^d} (T_t v)h\rd
\mu_* \right|=\sup_{h\in{\rm
L}^\infty\hskip-3pt,h\in[0,1]}\int_{{\mathbb R}^d}(T_t v)h
\rd \mu_* \\
&=\sup_{h\in{\rm L}^\infty\hskip-3pt,h\in[0,1]}\int_{{\mathbb
R}^d}(T_t h)v \rd \mu_* =\int_{{\mathbb R}^d}(T_t 1)v \rd \mu_*
=\int_{{\mathbb R}^d}v \rd \mu_* =\|v\|_1\mbox{.\qed}
\end{aligned}
\]

%%%%%%%%%%%%%%%%%%%%%%%%%%%%%%%%%%%%%%%%%%%%%%%%%%%%%%%%%%

\subsection{Linear case. Entropy Method}
\label{ssect.lcem}
The behaviour of the heat kernel of the linear
evolution considered and the L$^1$ contraction property allow to
notice that for all $t\ge t_0$
\[
\|u(t)\|_2^2\le \|u(t)\|_1\|u(t)\|_\infty\le C\frac{\|u_0\|_1}{t^{1/2}}.
\]
Notice that the above bound is sharp. In fact, consider the solution
corresponding to the Dirac delta at $x_0$, namely
$v(t,x)=K(t,x,x_0)$. Its L$^2$ norm then satisfies, using the
symmetry of the heat kernel and the semigroup property:
\[
\begin{split}
\Vert v(t)\Vert^2_2&=\int_{{\mathbb R}^d}K(t,x,x_0)^2\rd \mu_*
=\int_{{\mathbb R}^d}K(t,x,x_0)K(t,x_0,x)\rd \mu_*\\
&=K(2t,x_0,x_0)\sim c t^{-1/2}\ \ \ {\rm for\  large\ }t
\end{split}
\]

It is easy to get the same result by {\it entropy methods}. Although
this is not necessary in the present case due to the previous
calculations, this will serve as a model for the strategy of proof
used in the nonlinear setting, and will make already apparent the
role of the Nash inequalities proved before.
\begin{prop}\label{4.13}
Let $F(t)=\|v(t)\|_2^2$. Then $F(t)\le c t^{-1/2}$ for all $t>t_0$.
\end{prop}

\proof First consider nonnegative data. Having shown that the
L$^1$ norm of such solutions is conserved and, moreover, using the
fact that $I_{m_*}[v(t)]$ is decreasing as a function of $t$, we
get that $I_{m_*}[v(t)]/\|v(t)\|_1^2\le c_0$ for all positive $t$
and for some $c_0>0$. We are then allowed to use \eqref{GN} with
$r=2, s=1$, so that we have
\[
\frac{\rd F(t)}{\dt}=-I_{m_*}[v(t)] \le - c \frac{F^3}{\|v(t)\|^4_1}= - c
\frac{F^3}{\|v_0\|^4_1}
\]
Thus we get, integrating the above differential inequality:
\[
F(t)\le \tilde c\frac{\|v_{0}\|^2_1}{t^{1/2}}.
\]
The same decay holds true for all L$^1$ data, since we may write
$-(v_0)_-\le v_0\le (v_0)_+$ and use the order preserving property
of the evolution and the decay bound already proved for
nonnegative (or nonpositive) solutions. In fact, denoting by
$v_{\pm}(s)$ the time evolved of $(v_0)_\pm$, we have first that, by
comparison, $-v_-(s)\le v(s)\le v_+(s)$ and $v^2(s)= v^2_+(s)+v^2_-(s)$. This, together with
the above decay property for nonnegative solutions
$\|v_\pm(s)\|^2_2\le \tilde{c}\|(v_{0})_\pm\|_1 s^{-1/2}$, implies
that
\[\begin{split}
F(t)=\|v(t)\|^2_2
    &= \|v_-(t)\|^2_2+\|v_+(t)\|^2_2
     \le \frac{\tilde{c}\left(\|(v_{0})_-\|^2_1+\|(v_{0})_+\|^2_1\right)}{t^{1/2}}
     \le \tilde{c}\|v_0\|^2_1  t^{-1/2} \mbox{.\qed}
\end{split}
\]
%\noindent {\sl Proof.~}
%Assuming that $Q(f_0)$ and $M(f_0)$ are not zero, the right hand side takes the value infinity both as $t\to 0$ and $t\to\infty$ hence there is a minimum for one or several intermediate values of $t$.
%We want to take a particular value of $t$ that almost minimizes the above formula, and we want that value to correspond to the range of not small $t$ where $A(t)=c_2t^{-1/4}$.
%We assumed that $Q(f_0)/M(f_0)^2\le c_0$ for some $c_0>0$, then we consider the 1-parameter quantity
%\[
%t_\alpha=\alpha\left[\frac{M^2}{Q(f_0)}\right]^{2/3}
%\]
%and observe that, trivially,
%\[
%t_\alpha\ge 1\qquad\iff\qquad \alpha c_0^{2/3}\le 1.
%\]
%Plugging $t_\alpha$ into \eqref{Quasi.Nash}, gives then $\|f_0\|_2^2\le Z(t_\alpha)$ which
%can be rewritten in the form
%\begin{equation}
%\|f_0\|_2^2\le \left( \frac{2}{\alpha}+c_2\alpha^{1/2}\right)Q(f_0)^{1/3}M^{4/3}
%\end{equation}
%and this concludes the proof.\qed
%In the Lebesgue measure case, this inequality is called Nash inequality, and is equivalent to a whole family of GNI as we shall see in the sequel.
%%%%%%%%%%%%%%%%%%%%%%%%

%%%%%%%%%%%%%%%%%%%%%%%%%%%%%%%%%%%%%%%%%%%%%%%%%%%%%%%%%%%%%%%%%%%%%%%%%%%
\section{Nonlinear Entropy Method}\label{sect.nlem}

Once the linear flow has been examined and its behaviour
described, we  prepare the way for the proof of convergence with
rate of the nonlinear flow via a new version of the
entropy-entropy dissipation method. We shall use the entropy and
Fisher information introduced at the end of Section
\ref{sect.prel}. The results of this section hold for any $m<1$,
but the main interest is in employing them for the case $m=m_*$ as
is done in the subsequent section.  From now on we revert to the
notations for space, time and flow variables introduced in
sections \ref{sect.intro} and \ref{sect.prel}. Thus, $w=w(s,y)$.

\subsection{Comparing linear and nonlinear quantities. The Fisher information}
We have to prove the basic inequalities that relate the linear and
the nonlinear quantities of the entropy method. We start the
analysis by  a new inequality between linear and nonlinear Fisher
information, then we recall a Lemma of \cite{BBDGV} which compares
linear and nonlinear entropy. We shall write $V_*$ instead of
$V_{D_*}$. We put
\begin{equation}\label{Fisher}
\mathcal{I}_m[w]=
\int_{\RR^d}\left|\nabla\left(\frac{w^{m-1}-1}{m-1}V_*^{m-1}\right)\right|^2
V_* w \dy,
\end{equation}
which is the (nonlinear) Fisher information. It can be linearized, as done in \cite{BBDGV}, by letting $w=1+\varepsilon g V_*^{1-m}$
and taking the limit as $\varepsilon\to 0$. We obtain the linearized form of the
Fisher information, that takes the expression of the Dirichlet form typical
of the linearized equation
\begin{equation}
I_{m}[w]
    = \int_{\RR^d}\left|\nabla (w-1)V_*^{m-1}\right|^2 V_*\dy
    = \int_{\RR^d}\left|\nabla g\right|^2 V_* \dy
\end{equation}
the relation between $g$ and $w$ is $g=(w-1)V_*^{m-1}$; it is not restrictive to let
$\varepsilon=1$ in the sequel. \noindent The next Lemma compares
in a quantitative way $I_{m}$ and $\mathcal{I}_m$. This is a first
attempt that will be improved subsequently for $m=m_*$ and $m\ne
m_*$. We drop the subindex $m$ from both quantities for brevity.

\begin{lem}\label{Fisher-lin-nonlin}
Let $0<W_0\le w\le W_1<+\infty$, be a measurable function on $\RR^d$, with $W_0<1$ and $W_1>1$, and assume that $\mathcal{I}(w)<+\infty$. Then the following inequality holds true
\begin{equation}
I [w]
    \le k_1\mathcal{I}[w]
      + k_2\int_{\RR^d}g^4V_*^{4-3m}\dy\\
\end{equation}
for any $m<1$, where $g=(w-1)V_*^{m-1}$, $k_1= 2W_1^{3-2m}$ and $k_2$ depends only on $W_1, W_0$, $m$ and $d$.
%\[
%k_2(m,W_0,W_1)= \dfrac{4(2-m)}{2}W_1^{2(1-m)}\left[1 + \dfrac{9-m}{3}\overline{W}\right]
%        +\dfrac{1}{W_0^{1-m}\overline{W}}\left[\dfrac{1}{W_0} + \dfrac{1}{(1-m)\overline{W}}\right]
%\]
%with $\overline{W}=\min\big\{W_1-1,1-W_0\big\}/4$.
\end{lem}

\noindent {\sl Proof.~} We have $w-1=g V_*^{1-m}$. We
first re-write the Fisher information \eqref{Fisher} in the
following way:
\begin{equation}\begin{split}
\mathcal{I}[w]
    &= \int_{\RR^d}\left|\nabla\left(\frac{w^{m-1}-1}{m-1}V_*^{m-1}\right)\right|^2 V_* w \dy\\
        %&=  \int_{\RR^d}\left|\nabla\left(\frac{w^{m-1}-1}{(m-1)(w-1)}(w-1)V_*^{m-1}\right)\right|^2
        %            V_* w\dx\\
    &:=  \int_{\RR^d}\left|\nabla\left(A(w)(w-1)V_*^{m-1}\right)\right|^2 V_* w \dy,
\end{split}
\end{equation}
where we have defined
\[
A(w):=\frac{w^{m-1}-1}{(m-1)(w-1)}=\frac{a(w)}{w-1}.
\]
It is easy to check that $A(1)=1$, $A(w)>0$, and that $A(w)\to 0$,
when $w\to\infty$. Moreover,
\begin{equation}\label{A'}
A'(w)=\frac{w^{m-2}-A(w)}{w-1}\le 0
\end{equation}
since the function $a(w)=(w^{m-1}-1)/(m-1)$ is concave in $w$, so that its incremental quotient $A(w)$
(taken in $w=1$) is a non-increasing function of $w$. If we let $W_0\le w\le W_1$,
with $0<W_0\le 1$ and $1\le W_1< +\infty$, we then have the bounds
\begin{equation}\label{bounds.A}
W_1^{m-2}=a'(W_1)\le |A(w)| \le a'(W_0)=W_0^{m-2}.
\end{equation}
We shall also need estimates for $|A'(w)|$ for $W_0\le w\le W_1$
as above, and in fact it is easy to check that $A'(1)=(m-2)/2$  and that $A^\prime$ is
bounded away from zero.
%We have the following estimates
%\[
%|A'(w)| \le \left\{
%    \begin{array}{lll}
%        \dfrac{1}{W_0^{1-m}\varepsilon}\left[\dfrac{1}{W_0} + \dfrac{1}{(1-m)\varepsilon}\right],
%            &\mbox{if }w\in[W_0,1-\varepsilon)\subseteq [W_0,W_1]\\[4mm]
%        \dfrac{2-m}{2}\left[1 + \dfrac{9-m}{3}\varepsilon\right],
%            &\mbox{if }w\in[1-\varepsilon,1+\varepsilon])\subseteq [W_0,W_1]\\[2mm]
%        \dfrac{2}{(1-m)\varepsilon}\left[1 + \dfrac{1}{(1-m)\varepsilon}\right],
%            &\mbox{if }w\in(1+\varepsilon,W_1])\subseteq [W_0,W_1]\\
%    \end{array}
%\right.
%\]
%which can be summarized as:
%\begin{equation}\label{Est.A'}
%|A'(w)|
%    \le \dfrac{2-m}{2}\left[1 + \dfrac{9-m}{3}\varepsilon\right]
%        +\dfrac{1}{W_0^{1-m}\varepsilon}\left[\dfrac{1}{W_0} + \dfrac{1}{(1-m)\varepsilon}\right]
%    := k_0(m,W_0,W_1),
%\end{equation}
%since we can put $\varepsilon=\min\big\{W_1-1,1-W_0\}/4$.
Letting now $w$ be a function, noticing that $w-1=g
V_*^{1-m}$ and that \eqref{A'} can be rewritten as
$(w-1)A'(w)+A(w)=w^{m-2}$, we get
\begin{equation}\begin{split}
\nabla\left[A(w)(w-1)V_*^{m-1}\right]&=\nabla \big[A(w)g\big]
     = A(w)\nabla g + A'(w)\big[\nabla w \big] g \\
    &= A(w)\nabla g + A'(w)\big[\nabla (1+g V_*^{1-m}) \big] g\\
    &= A(w)\nabla g + A'(w)gV_*^{1-m}\big[\nabla g\big]
                    + A'(w)g^2 \big[\nabla V_*^{1-m} \big]\\
    &= \big[A(w)+ A'(w)(w-1)\big]\nabla g+ A'(w)g^2 \big[\nabla V_*^{1-m} \big]\\
    & = w^{m-2}\nabla g+ A'(w) \big[\nabla V_*^{1-m} \big]g^2.
\end{split}
\end{equation}
Now we can use this equality in \eqref{Fisher} to get:
\[
\begin{split}
\mathcal{I}[w]
    &= \int_{\RR^d}\left|A(w)\nabla g+\big[\nabla A(w)\big]g\right|^2 V_* w \dy\\
    &= \int_{\RR^d}\left|w^{m-2}\nabla g+ A'(w) \big[\nabla V_*^{1-m} \big]g^2 \right|^2 V_* w \dy\\
    &\ge \frac12 \int_{\RR^d}\left|\nabla g\right|^2 V_* w^{2(m-2)+1} \dy
        - \int_{\RR^d}g^4\big|A'(w)\big|^2 \left|\nabla V_*^{1-m}\right|^2 V_* w \dy\\
    &\ge \frac12W_1^{2m-3}\int_{\RR^d}\left|\nabla g\right|^2 V_* \dy
        - \int_{\RR^d}g^4\big|A'(w)\big|^2 \left|\nabla V_*^{1-m}\right|^2 V_* w \dy,
\end{split}
\]
where we have used the inequality $|a+b|^2 + |b|^2\ge (1/2)|a|^2
$ valid for any $a,b\in \RR$, and the bounds $W_0\le w\le W_1$.
Thus, we have
\[
I [g]=\int_{\RR^d}\left|\nabla g\right|^2 V_* \dy
    \le \frac{2}{W_1^{2m-3}}\mathcal{I}[w]
     + \frac{1}{W_1^{2m-3}}
        \int_{\RR^d}g^4\big|A'(w)\big|^2 \left|\nabla V_*^{1-m}\right|^2 V_* w \dy\\
\]
We next remark that the weight
\[\begin{split}
\left|\nabla V_*^{1-m}(y)\right|^2 V_*(y)
    &= \frac{4 |y|^2}{\big(D+|y|^2\big)^4} \frac{1}{\big(D+|y|^2\big)^{\frac{1}{1-m}}}\\
    &\le \frac{4}{\big(D+|y|^2\big)^3} \frac{1}{\big(D+|y|^2\big)^{\frac{1}{1-m}}}
    = \frac{4}{\big(D+|y|^2\big)^{3+\frac{1}{1-m}}}=4 V_*^{4-3m}
\end{split}
\]
is integrable whenever $(d-6)m>(d-8)$. \textit{Notice that when $m=m_*=(d-4)/(d-2)$, the weight is integrable.} We conclude by estimating $|A'|\le k_0$ so that
\[
I [g]=\int_{\RR^d}\left|\nabla g\right|^2 V_* \dy
    \le 2W_1^{3-2m}\mathcal{I}[w]
     + 4 W_1^{2(1-m)}k_0
        \int_{\RR^d}g^4V_*^{4-3m}\dy.\\
\]
This concludes the proof.\qed

\subsection{Evolution properties of the Fisher information}

We now describe some further properties of the Fisher information
$\mathcal{I}[w(s)]$ as a function of time, such as the fact that
it is uniformly bounded for large $s$ and it goes to zero as
$s\to+\infty$. We prove a new differential inequality for the
Fisher information. Indeed,  by Proposition 2.6 of \cite{BBDGV},
it is easy to see that the Fisher information  is finite almost
everywhere and is the time derivative of the entropy almost
everywhere so that:
\begin{equation}
\mathcal{F}[w(s_0)]-\mathcal{F}[w(s)]=\int_{s_0}^s\mathcal{I}[w(\xi)]\rd\xi
\end{equation}
taking the limits $s\to\infty$ and $s_0\to 0$, recalling that
$0\le \mathcal{F}[w(0)]<+\infty$, $0\le \mathcal{F}[w(s)]\to 0$ as
$s\to+\infty$, we can conclude that $\mathcal{I}[w(s)]$ is
integrable (and nonnegative) on $(0,+\infty)$.

%Since $w$ has $C^\infty$ regularity in space, it is easy to conclude, using the equation, that the first %time derivative is also continuous. Then also the Fisher information has a continuous time derivative, as %composition of continuous functions. We conclude that $\mathcal{I}[w(t)]$ is a nonnegative $C^1$ function %on $(0,\infty)$, which is integrable in time over $(0,+\infty)$. This implies that there is at least a %point $s_0\in(0,+\infty)$, so that $0\le \mathcal{I}[w(s_0)]<+\infty$. We can even choose $s_0$ to be the %smallest point $s_0\ge 0$ for which the Fisher information is finite. The next lemma shows that once the %Fisher information is finite at a point $s_0$, then it is finite for any $s\ge s_0$.
\begin{prop}\label{prop.diff.ineq.Fisher}
In addition to the running assumptions, suppose that $v(0)-V_D\ge0$.
Then, the following differential inequality for the Fisher
information holds true
\begin{equation}\label{diff.ineq.Fisher}
\frac{\rd \mathcal{I}[w(s)]}{\ds}\;\le\; \kappa_1 \mathcal{I}[w(s)]\;-\;\kappa_2 \mathcal{I}^2[w(s)]
\end{equation}
the constant $\kappa_1$ depends on $m,d,s_0$, and the constant $\kappa_2$ depends on $m,d$, the relative mass
$\int_{\RR^d}(v_0-V_D)\dx$, $W_0$ and $W_1$. Moreover, $\mathcal{I}[w(s)]$ goes to zero as $s\to\infty$.
\end{prop}

\noindent {\bf Remark.} The time derivative of the Fisher information is usually controlled by means of the Bakry-Emery method (cf. \cite{MR889476}) that allows to obtain spectral gap estimates. Such an estimate cannot hold when $m=m_*$ since there is no spectral gap. The above proposition can be viewed as a substitute for the Bakry-Emery method and gives a solution for asymptotic estimates in applications with no spectral gap.

\noindent {\sl Proof.~}The proof is divided in several steps. We use the notation
$\Omega=\left(\frac{v^{m-1}-V_D^{m-1}}{m-1}\right)$ in this section
for brevity. Note that for large $s$, $\Omega $ is uniformly bounded and \ $|\Omega|\le V_{D_0}^{m-2}|v-V_D|$.

\medskip

\noindent$\bullet~$\textsc{Expression of the derivative.} We first
perform a formal time-derivative of $\mathcal{I}$, but in this
case is convenient to write it in terms of $v$ and $V_D$ instead
of $w$, where we recall that $w=v/V_D$
\begin{equation}
\begin{split}
\frac{\rd \mathcal{I}}{\rd s}
    &=\frac{\rd}{\rd s}\int_{\RR^d}\left|\nabla\Omega\right|^2 v\dy\\
    &=  2\int_{\RR^d}\nabla\Omega
                \cdot\nabla\left(v^{m-2}\frac{\rd v}{\ds}\right)v\dy
        +\int_{\RR^d}\left|\nabla\Omega\right|^2 \frac{\rd v}{\ds}\dy
    = (A)+(B)
\end{split}
\end{equation}
Now we treat the two terms separately.

\noindent$\bullet~$\textsc{Estimating the term} (A): We have
\begin{equation}
(A) =2\int_{\RR^d} v\nabla\Omega
                \cdot\nabla\left(v^{m-2}\frac{\rd v}{\ds}\right)\dy
     =-2\int_{\RR^d} \nabla\cdot \left[v\nabla\Omega\right]
                v^{m-2}\frac{\rd v}{\ds}\dy.
\end{equation}
Using the equation, $v_s=\nabla\cdot \Omega$, we get
\begin{equation}
(A) =-2\int_{\RR^d}
        \left[\nabla\cdot \left(v\nabla\Omega\right)\right]^2
                v^{m-2}\dy
     =-2\int_{\RR^d}
        \left[\nabla\cdot \left(v\nabla\Omega\right)\right]^2
                \frac{\Omega^2}
                    {\Omega^2v^{2-m}}\dy.
\end{equation}
Then we have
\begin{equation}
\begin{split}
      (A)      & \le_{(i)} -2\frac {\left[\int_{\RR^d}
            |\nabla\cdot \big( v\nabla\Omega\big)|
            |\Omega|\dy \right]^2}
            {\int_{\RR^d}\Omega^2v^{2-m}\dy}
           \le_{(ii)}-2\frac{\left[-\int_{\RR^d}
            v|\nabla\Omega|^2 \dy \right]^2}
                    {\int_{\RR^d}\Omega^2v^{2-m}\dy}\\
                   &\le_{(iii)}-2\frac{\mathcal{I}^2}{c_2\int_{\RR^d}(v(0)-V_D)\dy}
                :=-\kappa_2\mathcal{I}^2,
\end{split}
\end{equation}
where in (i) we have used the H\"older inequality
\[
\int \frac{h_1^2}{h_2}\rd\mu\ge\frac{\left[\int h_1\rd\mu\right]^2}{\int h_2\rd\mu},
\]
while in (ii) we use  integration by parts, after noticing that
$|a||b|\ge ab$. The point (iii) relies on fact that the difference
between two Barenblatt solutions behaves like $V_D^{2-m}$, and on
the fact that $V_{D_0}\le v(t)\le V_{D_1}$, so that
\[
\begin{split}
\int_{\RR^d}\Omega^2v^{2-m}\dy
    &=\int_{\RR^d}\left(\dfrac{w^{m-1}-1}{m-1}\right)^2V_D^{2(m-1)}v^{2-m}\dy\\
    &=\int_{\RR^d}\left(\dfrac{w^{m-1}-1}{m-1}\right)^2V_{D_1}^{m}\dy\\
 (a)&\le\max\{W_0^{m-2},W_1^{m-2}\}\int_{\RR^d}|w-1|^2V_{D_1}^{m}\dy\\
    &=\max\{W_0^{m-2},W_1^{m-2}\}\int_{\RR^d}|v-V_D|^2\frac{V_{D_1}^{m}}{V_D^{2}}\dy\\
 (b)&\le c_0\max\{W_0^{m-2},W_1^{m-2}\}\int_{\RR^d}|v-V_D|^2V_{D_1}^{m-2}\dy\\
 (c)&\le c_0c_1\max\{W_0^{m-2},W_1^{m-2}\}\int_{\RR^d}|v-V_D|V_{D}^{2-m}V_{D}^{m-2}\dy\\
 (d)& = c_2 \int_{\RR^d}|v-V_D|\dy = c_2 \int_{\RR^d}(v-V_D)\dy =
  c_2\int_{\RR^d}(v(0)-V_D)\dy,
\end{split}
\]
where in (a) we have used \eqref{bounds.A}, namely
\[
W_1^{m-2}|w-1|\le \left|\frac{(w^{m-1}-1)}{(m-1)}\right| \le
W_0^{m-2}|w-1|,
\]
while in (b) we have used $V_D\ge c_0{V_{D_1}}$ and in (c) we have
used $|v-V_D|\le c_1 V_D^{2-m}$. In the last step (d) we have used
hypothesis (H2') together with the fact that $v(0)-V_D\ge 0$ and
conservation of relative mass, proved in Proposition 2.3 of
\cite{BBDGV}.

\noindent$\bullet~$\textsc{Estimating the term} (B). We shall use the celebrated B\'enilan-Crandall estimates \cite{Benilan-Crandall}, that for solutions to the un-rescaled FDE $\partial_t u =\Delta u^m/m$ read
\[
\partial_t u(t,x) \le \frac{u(t,x)}{(1-m)t}\qquad\mbox{for any $t>0$}
\]
if $m<1$, even for $m\le 0$. We perform the scaling to the
Fokker-Plank equation, like in section \ref{sect.prel}, so that
the B\'enilan--Crandall estimates read
\begin{equation}\label{BC.est.FP}\begin{split}
\partial_s v(s,y)
    &\le\frac{2}{[d(1-m)-2](1-m)}\left[\frac{d}{d(1-m)-2}
            +\frac{1}{(1-m)\left(\ee^{s(1-m)[d(1-m)-2]/2}-1\right)}\right]v(s,y)\\
    &\le\frac{2}{[d(1-m)-2](1-m)}\left[\frac{d}{d(1-m)-2}
            +\frac{1}{(1-m)\left(\ee^{s_0(1-m)[d(1-m)-2]/2}-1\right)}\right]v(s,y)\\
    &=\kappa_1(m,d,s_0)v(s,y),
\end{split}
\end{equation}
if $s\ge s_0>0$. We remark that $\kappa_1\to +\infty$ when $s_0\to
0$ but this will not be a problem. We finally estimate (B)
\begin{equation}
\begin{split}
(B)  &=\int_{\RR^d}\left|\nabla\Omega\right|^2 \partial_s v\dy
      \le \mathcal{C}(m,d,s_0)
        \int_{\RR^d}\left|\nabla\Omega\right|^2 v\dy\\
     &=\;\kappa_1(m,d,s_0)\;\mathcal{I}[w(s)].
\end{split}
\end{equation}
This calculation is formal and has to be justified, but  before we
do that let us draw a first consequence.

\medskip

\noindent$\bullet~$\textsc{Integrating the Differential
Inequality.} We obtained a closed differential inequality for the
Fisher information $\mathcal{I}[w(s)]=\mathcal{I}(s)$
\[
\frac{\rd\mathcal{I}(s)}{\ds}-\kappa_1\mathcal{I}(s)+\kappa_2\mathcal{I}^2(s)\le 0
\]
which is of Bernoulli type and can be estimated explicitly.
Indeed, the exact solution on $(s_1,s)\subseteq [0,+\infty)$ of
the Bernoulli ordinary differential equation
$Z^\prime(s)-\kappa_1Z(s)+\kappa_2Z^2(s)=0$ is given by
\begin{equation}
Z(s)=\frac{\ee^{\kappa_1(s-s_1)}}{\left[Z_0^{-1}
                +\int_{s_1}^s\ee^{\kappa_1(\xi-s_1)}\kappa_2\rd\xi\right]}
            \le\frac{\ee^{\kappa_1(s-s_1)}}
                {\kappa_2\int_{s_1}^s\ee^{\kappa_1(\xi-s_1)}\rd\xi}\le c\frac{\kappa_1}{\kappa_2}
\end{equation}
if $s\ge s_1+1:=s_0$, for a suitable $c>1$ which can be taken to
be arbitrarily close to one by choosing $s_1$ large enough. By
comparison it is clear that $\mathcal{I}(s)\le Z(s)\le
c\kappa_1/\kappa_2$, provided $\mathcal{I}(s_0)\le Z(s_0)$.
Therefore, for all $c>1$ and $0<s_0\le s$,
\begin{equation}
\mathcal{I}[w(s)]\le\frac{c\kappa_1}{\kappa_2}.
\end{equation}
The constant $\kappa_2$ depends on $m,d$, the relative mass
$\int_{\RR^d}(v_0-V_D)\dx$, $W_0$ and $W_1$; the constant $\kappa_1$
depends on $m,d,s_0$ and $\kappa_1\to+\infty$ when $s_0\to 0$.

\medskip

\noindent$\bullet~$\textsc{Justification of the calculation.} The
differentiation of $\mathcal{I}$ performed above contains
calculations that are not justified in principle since they
involve differentiations and  integrations by parts in integrals
over the whole space that are not justified a priori. Therefore,
we introduce a cutoff function $\zeta_n(y)$ for the integrand of
$\mathcal{I}$ and define
$$
\mathcal{I}_n=\int_{\RR^d}\left|\nabla\Omega\right|^2 v
\zeta^2_n\dy.
 $$
We assume that $\zeta_n$ has  value $1$ if $|y|\le n$, value 0 of
$|y|\ge 2n$, and $|\nabla \zeta_n|\le 1/n$, $ |\Delta\zeta_n|\le
1/n^2$. Then we have
\begin{equation}
\frac{\rd \mathcal{I}_n}{\rd s}= (A_n)+(B_n)
\end{equation}
and the two terms are as before but for the cutoff factor. There
is no problem with $(B_n)$. But $(A_n)$ produces extra terms that
we must control. Indeed,
\begin{equation*}
\begin{split}
(A_n) &=2\int_{\RR^d} v\nabla\Omega
                \cdot\nabla\left(v^{m-2}\frac{\rd
                v}{\ds}\right)\zeta_n^2\dy\\
     &=-2\int_{\RR^d} \nabla\cdot \left(v\nabla\Omega\right)
                v^{m-2}\frac{\rd v}{\ds}\zeta_n^2\dy
     -2\int_{\RR^d}
                v^{m-1}\frac{\rd v}{\ds}\,\left(\nabla\Omega\cdot \nabla \zeta_n^2\right)\dy.
     \end{split}
\end{equation*}
When we replace $dv/ds$ by its value according to the equation,  the
first of the two terms of the last expression becomes
\begin{equation}\label{form.An1}
(A_{n1}) :=-2\int_{\RR^d}
        \left|\nabla\cdot \left(v\nabla\Omega\right)\right|^2
                v^{m-2}\zeta_n^2\dy \left(=-2\int_{\RR^d}|v_s|^2 v^{m-2}\zeta_n^2\dy\right).
\end{equation}
which has a convenient negative sign. We
now perform integration by parts on this term, an operation that is
now perfectly justified, and we get much as before:
\begin{equation*}
\begin{split}
      (A_{n1}) &=-2\int_{\RR^d}
        \left[\nabla\cdot \left(v\nabla\Omega\right)\right]^2
                \frac{\Omega^2}
                    {\Omega^2v^{2-m}} \zeta_n^2\dy
\le -2\frac{\left[\int_{\RR^d}
            \left|\nabla\cdot \left(v\nabla\Omega\right)\right|
            |\Omega|\zeta_n^2\dy \right]^2}
                    {\int_{\RR^d}\Omega^2v^{2-m} \zeta_n^2\dy}
                    \end{split}
\end{equation*}
The numerator of the last term is larger than $|\int_{\RR^d}
\nabla\cdot \left(v\nabla\Omega\right) \Omega\zeta_n^2\dy|$, hence
\begin{equation*}
(A_{n1})\le -\kappa_2 \left| \int_{\RR^d} \left(\nabla\cdot
\left(v\nabla\Omega\right)\right)
                           \Omega\zeta_n^2\dy\right|^2.
\end{equation*}
with the notation that we have used above. Let us calculate the
integral: after integrating by parts, it gives a term as before plus
a term where $\zeta_n^2$ is differentiated, as follows:
\begin{equation*}
\int_{\RR^d} v\left|\nabla\Omega\right|^2\zeta_n^2 \dy +2
\int_{\RR^d} v\Omega\nabla\Omega
\cdot\zeta_n\nabla\zeta_n\dy=(X_n')+(X_n'')
\end{equation*}
The first term is $(X_n')=\mathcal{I}_n$, as before, while the new
term,  $(X_n'')$, can be tackled as follows. We separate by H\"older
a factor like $\mathcal{I}_n^{1/2}$ (but we only need to integrate
in the annulus $R_n=\{n\le |y|\le 2n\}$ so it goes to zero as
$n\to\infty$) and we still have another factor:
\begin{equation*}
\int_{\RR^d} v\left|\Omega\right|^2 |\nabla\zeta_n|^2\dy\le C
\int_{R_n}v |V_{D_0}^{m-1}-V_{D_1}^{m-1}|^2 |\nabla\zeta_n|^2\dy\le
C\int_{R_n} \frac{v}{n^2}dy
\end{equation*}
and this tends to zero as $n\to\infty$ for $m>m_*$. For $m\le m_*$
we calculate differently,
\begin{equation*}
\int_{R_n}v |V_{D_0}^{m-1}-V_{D_1}^{m-1}|^2 |\nabla\zeta_n|^2\dy\le
C\int_{R_n} \frac{V_D^{m-1}}{n^2}|v-V_D|dy\le C \int_{R_n} |v-V_D|dy,
\end{equation*}
that goes to zero as $n\to\infty$, but in a uniform way we only know
that is bounded a priori. In any case, raising to the square we get
an estimate of the form
\begin{equation}
(A_{n1})\le -\kappa'_2 \mathcal{I}_n^2+ \kappa_2''\mathcal{I},
\end{equation}
with constants uniform in $n$.

\medskip

\noindent $\bullet$ We now consider the other new
term
\begin{equation*}
      (A_{n2}) =-4\int_{\RR^d} v^{m-1}\zeta_n
           \frac{\rd v}{\ds}  \left( \nabla\Omega \cdot \nabla \zeta_n \right)\, \dy.
\end{equation*}
Use H\"older and the numerical inequality $2ab\le \varepsilon a^2
+ b^2/\varepsilon$, to separate a term like $(A_{n1})$ in formula
\eqref{form.An1}, with a factor $\varepsilon$ which is convenient
to be absorbed by the term already existing that has a negative
sign in front, as we have remarked. We are then left with another
factor of the form
\begin{equation*}
   (A_{n22})   \le   C\int_{\RR^d} v^{m}\left|\nabla\Omega\right|^2
                | \nabla \zeta_n|^2\dy\sim C \int_{n\le |y|\le 2n}\frac{v^m}{n^2}
     |\nabla (v^{m-1}-V^{m-1})|^2\dy.
\end{equation*}
Since $v^{m-1}(s,y)\sim |y|^2$ as $|y|\to\infty$, and we get an
equivalent expression
\begin{equation*}\label{formi2n}
   (A_{n22})   \le   C\int_{n\le |y|\le 2n} \left|\nabla\Omega\right|^2
                v\,\dy\le C \mathcal{I}_{2n}.
\end{equation*}
If we had $C \mathcal{I}_{n}$ instead of $C \mathcal{I}_{2n}$  we
would have ended. The estimate for $(B_n)$ has no problems.

\medskip

\noindent $\bullet$ To solve the difficulty we take a little
detour. We integrate the inequality obtained so far for
$\rd\mathcal{I}_{n}/\ds$ to get the integrated inequality:
\begin{equation}
 \mathcal{I}_{n}(s_2)- \mathcal{I}_{n}(s_1)\le k \int_{s_1}^{s_2} \mathcal{I}_{n}\ds
 + C\int_{s_1}^{s_2} \mathcal{I}_{2n}\ds
\end{equation}
But the right-hand side is bounded above by the integral
\begin{equation}
(C+k)\int_{s_1}^{s_2} \mathcal{I}\,\ds
\end{equation}
and this is known to be bounded  by the relative entropy.
Moreover, since the integral $\int_{s_1}^{\infty}
\mathcal{I}\,\ds$ is finite, for every $\ve>0$ there exists a
$s_\ve$ such that
\begin{equation}\label{integralofI}
\int_{s_\ve}^{\infty} \mathcal{I}\,\ds\le \ve.
\end{equation}
We conclude that $\mathcal{I}_{n}(s_2)- \mathcal{I}_{n}(s_1)\le
\ve$ when $s_\ve\le s_1\le s_2$. Combining this half continuity
with the integrability of $\mathcal{I}_n(s)\le \mathcal{I}(s)$
given by \eqref{integralofI}, we obtain by an easy calculus lemma
that
\begin{equation}
\mathcal{I}_n(s)\le C_1\ve
\end{equation}
for all $s\ge 2 s_\ve$, with $C_1$ uniform in $n$. We conclude
that $\mathcal{I}(s)=\lim_n \mathcal{I}_n(s)$ is bounded for all
large times and goes to zero as $s\to\infty$.

\noindent $\bullet$ Coming back to the differential inequality
satisfied by ${\mathcal I}_n$ we have proved that ${\mathcal
I}_n^\prime\le c_1{\mathcal I}-c_2 {\mathcal I}_n^2$. Integrating
this differential inequality in time between $s_1$ and $s_2$ with
$s_1<s_2$ sufficiently large we get $I_n(s_2)-I_n(s_1)\le
c_1\int_{s_1}^{s_2}I(s)\ds-c_2\int_{s_1}^{s_2}I_n(s)^2\ds$ so
that, passing to the limit as $n\to+\infty$ and using both
monotone convergence and the boundedness of $I$ as a function of
time we get $I(s_2)-I(s_1)\le
c_1\int_{s_1}^{s_2}I(s)\ds-c_2\int_{s_1}^{s_2}I(s)^2\ds$, which is
an equivalent form of our statement.\qed

%%%%%%%%%%%%%%%%%%%%%%%%%%%%%%%%%%%%%%%%%%%%%%%%%%%%%%%%%%%%%%%%%%%%%%%%%%%%%%
\subsection{Comparing linear and nonlinear entropies}
The quantitative comparison of  linear and nonlinear entropies
concludes the preliminary results needed for the nonlinear entropy
method. Under Assumptions (H1'')-(H2''), the relative entropy is
well defined.

\begin{lem}[An equivalence result]\label{Lem.Bounds.RE}
Let $m<1$. If $w$ satisfies {\rm (H1'')-(H2'')}, then
\begin{equation}\label{Disug.Entr.Lin-Nolin}
\frac{F[w]}{2 W_1^{2-m}} \le \mathcal{F}[w]\le \frac{F[w]}{2
W_0^{2-m}}.
\end{equation}
We recall that $F[w]= \int_{\RR^d}|w-1|^2V_{D_*}^{m}\dx $.
\end{lem}

 The
short proof of this result has been given first in \cite{BBDGV}
but we repeat it here for reader's convenience. \noindent\proof
For $a>0$, let
$\phi_{a}(w):=\frac{1}{1-m}\left[(w-1)-(w^m-1)/m\right]
-a\left(w-1\right)^2$. We compute
$\phi_{a}'(w)=\frac{1}{1-m}\left[1-w^{m-1}\right]-2a\left(w-1\right)$
and $\phi_{a}''(w)=w^{m-2}-2a$, and note that
$\phi_{a}(1)=\phi_{a}'(1)=0$. With $a=W_1^{m-2}/2$, $\phi_{a}''$
is positive on $(W_0,W_1)$, which proves the lower bound after
multiplying by $V_D^m$ and integrating over $\RR^d$. With
$a=W_0^{m-2}/2$, $\phi_{a}''$ is negative on $(W_0,W_1)$ which
proves the upper bound. \qed

Equivalently, we may write
\begin{equation}\label{Disug.Entr.Lin-Nolin.time}
\frac{F[w]}{2 W_1^{2-m}}\le \frac{F[w]}{2\sup\limits_{\RR^d}|w|^{2-m}}
\le \mathcal{F}[w]\le \frac{F[w]}{2 \inf\limits_{\RR^d}|w|^{2-m}}\le \frac{F[w]}{2 W_0^{2-m}}.
\end{equation}

\subsection{The entropy bounds a suitable $\LL^p$-norm}

\begin{lem}\label{Lem.Entr.Lp}
Let $m<1$. If $w$ satisfies {\rm (H1'')-(H2'')}, then
\begin{equation}\label{Disug.Entr.Lp}
\|w-1\|^{2+\frac{m}{1-m}}_{\LL^{2+\frac{m}{1-m}}(\RR^d)}
    \le \overline{D}_m F[w],
\end{equation}
where $\overline{D}_m$ is given at the end of the proof.
\end{lem}

\noindent {\sl Proof.~}
We first state some inequalities between Barenblatt solutions with different constants.
Consider
\[
\frac{\partial V_D}{\partial D}
    =-\frac 1{1-m}\left[D+|y|^2 \right]^{-\frac{2-m}{1-m}}
    =-\frac 1{1-m}V_D^{2-m}\le 0\;.
\]
Hence, for any $0<D_1<D_0$
\begin{equation*}\label{est.diff.baren}
\frac{D_0-D_1}{1-m}V_{D_0}^{2-m}\le\left|V_{D_1}-V_{D_0}\right|
    \le \frac{D_0-D_1}{1-m}V_{D_1}^{2-m}
\end{equation*}
Moreover, it is easy to see that if $0<D_1\le D_0$
\[
V_{D_0}^{1-m}(y)=\frac{1}{D_0+|y|^2}
    \le \frac{1}{D_1+|y|^2}=V_{D_1}^{1-m}(y)
    \le \left(1+\frac{D_0}{D_1}\right)\frac{1}{D_0+|y|^2}
    =\left(1+\frac{D_0}{D_1}\right)V_{D_0}^{1-m}(y).
\]
The above inequalities prove that
$|w-1|^{m/(1-m)}$ is bounded by a multiple of
$V_*^m$. Indeed, by hypothesis (H1') we have that $0<D_0<D_*<D_1$,
so that $V_{D_0}-V_{*}\le v(s)-V_{*}\le V_{D_1}-V_{*}$, and
since $w=v/V_{*}$
\[\begin{split}
|w-1|
    &=\left|\frac{v(s)-V_{*}}{V_{*}}\right|
     \le\frac{|V_{D_1}-V_{D_0}|}{V_{*}}
     \le \frac{D_0-D_1}{1-m}         \frac{V_{D_1}^{2-m}}{V_{*}}
\end{split}
\]
Thus,
\[\begin{split}
\|w-1\|^{2+\frac{m}{1-m}}_{\LL^{2+\frac{m}{1-m}}(\RR^d)}
    &=\int_{\RR^d}|w-1|^2|w-1|^{\frac{m}{1-m}}\dy\\
    &\le\left(\frac{D_0-D_1}{1-m}\right)^{\frac{m}{1-m}}
         \left(1+\frac{D_*}{D_1}\right)^{\frac{m(2-m)}{(1-m)^2}}
         \int_{\RR^d}|w-1|^2V_{*}^{m}\dy:= \overline{D}_m F[w]
         \mbox{.\qed}
\end{split}
\]
\noindent\textbf{Remarks. }  (i) The estimate proves that
$w-1\in\LL^{2+\frac{m}{1-m}}(\RR^d)$, whenever the initial entropy
is finite, since we know, joining inequalities
\eqref{Disug.Entr.Lin-Nolin} and \eqref{Disug.Entr.Lp}:
\begin{equation}\label{Lp-entr}
\|w(s)-1\|^{2+\frac{m}{1-m}}_{\LL^{2+\frac{m}{1-m}}(\RR^d)}
    \le F[w(s)]
    \le 2\overline{D}_mW_1^{2-m}\mathcal{F}[w(s)]
    \le 2\overline{D}_mW_1^{2-m}\mathcal{F}[w_0],
\end{equation}
since the nonlinear entropy is decreasing in time. Moreover, we
have also proved that
\begin{equation}
\|w(s)-1\|^{2+\frac{m}{1-m}}_{\LL^{2+\frac{m}{1-m}}(\RR^d)}
    \le 2\overline{D}_mW_1^{2-m}\mathcal{F}[w(s)]
\end{equation}
and  we shall show below that entropy goes to zero as $s\to+\infty$\,.

\noindent (ii) As an easy consequence, letting $w-1=(v-V_*)/V_*$
and using the fact that $V_*\le C$, we obtain
\begin{equation}\label{Disug.Entr.Lq}
\|v-V_*\|^{2+\frac{m}{1-m}}_{\LL^{2+\frac{m}{1-m}}(\RR^d)}
    \le \overline{D}_m F[w]\,.
\end{equation}
(iii) For $m=m_*$ we have $2+\frac{m}{1-m}= d/2$.

\subsection{Comparing linear and nonlinear Fisher information}

With the above remarks we can improve on Lemma
\ref{Fisher-lin-nonlin} that compares the linear and nonlinear
Fisher information:

\begin{prop}\label{Fisher-lin-nonlin.2}
Under the same assumptions of Lemma {\rm \ref{Fisher-lin-nonlin},}
we have
\begin{equation}\label{diseg.Fisher-lin-nonlin.2}
I[g]=\int_{\RR^d}\left|\nabla g\right|^2 V_* \dy
    \le k_1\mathcal{I}[w]
      + k_3\mathcal{F}^{1+\sigma}[w]\\
\end{equation}
for any $m<1$, where $g=(w-1)V_*^{m-1}$, $\sigma=2/[d+2+m/(1-m)]>0$,  $k_1= 2W_1^{3-2m}$, and
$k_3>0$ is given at the end of the proof.
\end{prop}
\noindent {\sl Proof.~} We estimate the second term of the inequality of Lemma \ref{Fisher-lin-nonlin} in the following way
\[\begin{split}
\int_{\RR^d}g^4V_*^{4-3m}\dy
    &= \int_{\RR^d}\big(|w-1|V^{m-1}\big)^4V_*^{4-3m}\dy
     = \int_{\RR^d}\big(|w-1|\big)^4V_*^{m}\dy\\
    &\le \left\|w-1\right\|_{\infty}^2\int_{\RR^d}\big(|w-1|\big)^2V_*^{m}\dy
     = \left\|w-1\right\|_{\infty}^2F[w]\\
\end{split}
\]
Now we recall the interpolation inequality \eqref{interp.Cj.Lp} with $j=0$
\begin{equation}
\|f\|_{\LL^\infty(\RR^d)}\;
\le\;\mathcal{C}_{d}\;\|f\|_{C^{1}(\RR^d)}^{\frac{d}{d+p}}\;\|f\|_p^{\frac{p}{d+p}}
\end{equation}
then we apply it to $f=w-1$ and we let $p=2+m/(1-m)$. We get:
\begin{equation}\label{infinity-interp}
\begin{split}
\|w-1\|_{\LL^\infty(\RR^d)}\;
    &\le\mathcal{C}_{p,d}\;\|w-1\|_{C^{1}(\RR^d)}^{\frac{d}{d+p}}
        \;\left(\|w(s)-1\|^{2+\frac{m}{1-m}}_{\LL^{2+\frac{m}{1-m}}(\RR^d)} \right)^{\frac{1}{d+p}}\\
    &\le\mathcal{C}_{p,d}M_1\left(2\overline{D}_mW_1^{2-m}
                \mathcal{F}[w]\right)^\frac{1}{d+2+\frac{m}{1-m}}:=k_3\mathcal{F}[w]^{\sigma/2}\\
\end{split}
\end{equation}
where $\sigma=2/[d+2+m/(1-m)]>0$ for any $m<1$ and we used inequality \eqref{Lp-entr} and the fact that $\|w-1\|_{C^{1}(\RR^d)}^{\frac{d}{d+p}}\le M_1$ by Theorem \ref{lem:holdereg}.
Thus we have proved that
\[\begin{split}
\int_{\RR^d}g^4V_*^{4-3m}\dy
    &\le\left\|w-1\right\|_{\infty}^2F[w]
        \le k_3\mathcal{F}[w]^{1+\sigma}.
\end{split}
\]
The expression of $k_3$ is then
\[
k_3=\mathcal{C}_{p,d}M_1\left(2\overline{D}_mW_1^{2-m}\right)^\frac{1}{d+2+\frac{m}{1-m}}
\]
where $\|w-1\|_{C^{1}(\RR^d)}^{\frac{d}{d+p}}\le M_1$,
$\overline{D}_m = \frac{D_0-D_*}{1-m}\frac{D_*-D_1}{1-m} \left(1+\frac{D_*}{D_1}\right)^{\frac{2-m}{1-m}}$, and $\mathcal{C}_{p,d}$ is the constant of the interpolation inequality \ref{interp.Cj.Lp} with $j=0$ and $p=2+m/(1-m)$.
\qed

\noindent\textbf{Remarks. } (i) The above proposition holds for any
$m<1$ and allow to conclude that $I(s)\to 0$ as $s\to +\infty$,
since we already know that both $\mathcal{I}(s)$ (cf. Proposition
\ref{prop.diff.ineq.Fisher}) and $\mathcal{F}(s)$ tend to zero as
$t\to +\infty$.

(ii) When $m=m_*$, we obtain that $\sigma=4/(3d)$. But in this critical case we shall need another finer
comparison between the linear and the nonlinear Fisher information that hold only when $m=m_*$, namely we would like to have that there exists $s_0>0$ and a constant $k_4>0$, such that
\begin{equation}\label{I.3}
    I[g(s)]\le k_4\mathcal{I}[w(s)].
    \end{equation}
for any $s\ge s_0$, where $g=(w-1)V_*^{m_*-1}$. Unfortunately the
above inequality is not guaranteed for all times $s\ge s_0$. In the
next section we will prove a weaker version of this statement,
sufficient to our scopes, namely we will show that the above
estimate \eqref{I.3} holds on a family of intervals $[s_{1,k},
s_{2,k}]$ that is sufficiently dense as $s\to\infty$.
The technical details will be postponed to Appendix A4.

%%%%%%%%%%%%%%%%%%%%%%%%%%%%%%%%%%%%%%%%%%%%%%%%%%%%%%%%%%%%%%%%%

\section{Proofs of the main results in the critical case}
\label{sect.nlem2}

In this section we shall always take $m=m_*$, and we shall show
that the nonlinear flow converges with the same rate as the linear
case, cf. Section \ref{ssect.lcem}. We shall use the relationship
between the entropy functional $\mathcal{F}$ and the Fisher
information $\mathcal{I}$, namely
$\rd\mathcal{F}/{\ds}=-\mathcal{I}$. In view of the absence of any
spectral gap (or Hardy-Poincar\'e inequality) inequality, valid
instead in the case $m\neq m_*$, we have to proceed differently.
The Gagliardo-Nirenberg inequality, that in the linear case give
the correct decay of the linearized entropy in the
$\LL^2(V_*^{2-m}\dx)$-norm, turn out to work as well in the
nonlinear case as the previous proposition started to show.

\medskip

\noindent{\bf Proof of Theorem \ref{thm.conv.entropy}.} Notice first
that $\|g\|_\infty$ is finite and bounded as a function of time. In
fact, by hypothesis (H1') we know that
\[\begin{split}
    |g(s,y)|&=|w(s,y)-1|V_{D_*}^{m-1}=\left|\frac{v-V_{D_*}}{V_{D_*}}\right|V_{D_*}^{m-1}
        \le c_0 |V_{D_1}(y)-V_{D_0}(y)|V_{D_*}^{m-2}\\
        &\le c_1 V_{D_*}^{2-m}V_{D_*}^{m-2} = c_1
\end{split}\]
for all $y\in\RR^d$ and all $s>0$, where $c_i$ are a positive
constant depending only on $m, D_0, D_1, D_{*}$. The inequality
$|V_{D_1}-V_{D_0}|\le c V_{D_*}^{2-m}$ can be proved easily using the
explicit expression of the pseudo--Barenblatt solutions (see the
proof of Lemma \ref{Lem.Entr.Lp}).

\medskip

Next we prove that $I$ is bounded as a function of time.  Indeed
by Lemma \ref{Fisher-lin-nonlin} we observe that
\begin{equation}\label{I.1}
    I[g]\le k_1\mathcal{I}[w]
          + k_2\int_{\RR^d}g^4V_{D_*}^{4-3m}\dy\\
        \le k_1\mathcal{I}[w]
          + k_2k_3\|g\|_\infty^4
\end{equation}
where we have noticed that, for $m=m_*$, $V_{D_*}^{4-3m}=\big(1+|x|^2\big)^{-(d+4)/2}$ is integrable. It has been proved in Proposition \ref{prop.diff.ineq.Fisher} that $\mathcal{I}$ is also bounded.

By conservation of relative mass, cf. Proposition 2.3 of
\cite{BBDGV}, we know that $\|g(s)\|_1=\|g(0)\|_1$, where we have
used the fact that $v_0-V_{D_*}$ above is taken also nonnegative with
$\int (v_0-V_{D_*})\,\dy=M>0$ in this part of the proof. This implies
that the ratio $I/M=I_{m_*}[g(s)]/\|g(s)\|_1^2$ is bounded as a
function of time.

We shall use now the Gagliardo-Nirenberg inequalities of
Proposition \ref{GNprop} taking $v=g(t)$, putting
$F=\|g\|^2_{L^2(V^{2-m}\dx)}$, $I$ the linear Dirichlet form, and
$M=\|g\|_{L^1(V^{2-m}\dx)}$:
\begin{equation}\label{GNI.3}
        F^3\le K_1 IM^4.
\end{equation}
The validity of such inequalities
depends on the boundedness of the ratio
$I_{m_*}[g(s)]/\|g(s)\|_{L^1(V^{2-m}\dx)}^2$, which is ensured
along the evolution, as above mentioned.

 We now prove  an entropy - entropy production
inequality. We obtain a differential inequality for the entropy
$\mathcal{F}$, by comparing it with the Fisher information
$\mathcal{I}$ via Gagliardo-Nirenberg inequalities,
\begin{equation}\label{entr.prod}
\begin{split}
\mathcal{F}^3[w]
    &\le^{(a)} \left[\frac 12W_0^{m-2}\int_{\RR^d}|w-1|^2V_{D_*}^{m}\dy\right]^3
       = \left[\frac 12W_0^{m-2}\right]^3F^3\\
    &\le^{(b)} \left[\frac 12W_0^{m-2}\right]^3K_1 IM^4=K_2IM^4\\
\end{split}
\end{equation}
where (a) follows from \eqref{Disug.Entr.Lin-Nolin} of Lemma
\ref{Lem.Bounds.RE}, while in (b) we used the Gagliardo-Nirenberg
inequality \eqref{GNI.3} above.

\medskip

\noindent (i) In order to  continue the argument, we
assume for the moment that the initial datum satisfies
$v_0\ge V_{D_*}$ and is radially symmetric so that
$g_0=(w_0-1)V_{D_*}^{m-1}$ is {\it nonnegative}. This extra assumption will be removed
afterwards. Under it we will prove in Appendix A4 that there is
an infinite sequence of intervals of  times $[s_{1,k}\,,
s_{2,k}]\subset[2k,2k+2]$ (hence, $s_{2,k}\le s_{1,k+1}$) such that
\begin{equation}
I[g(s)]\le k_4\,\mathcal{I}[g(s)]\qquad\mbox{\rm for all} \ s\in
\bigcup_{k\in\NN}[s_{1,k}\,, s_{2,k}]
\end{equation}
for a constant $k_4$ that does not change along the evolution. We shall
prove moreover that the length of each of such intervals is at least
$1/2$ for all $k\ge k_0$, which in particular implies that
\begin{equation} \sum_{k=k_0}^n
(s_{2,k}- s_{1,k})\ge \sum_{k=k_0}^n \frac{1}{2}=\frac{n-k_0}{2}\ge
c\,\normalcolor s_{2,n}
\end{equation}
whenever $n\ge n_0$ is large, for a suitable $c>0$. Then, recalling that $\mathcal{I}=-\rd\mathcal{F}/\ds$, and
using \eqref{entr.prod} we conclude that
\[
\mathcal{F}^3\le K_2\,k_4\,M^4\, I \le k_5 \mathcal{I}= -k_5\frac{\rd\mathcal{F}}{\dt}\,,
\]
and an integration oven the interval $[s_{1,k}, s_{2,k}]$
gives
\[
\sum_{k=1}^{n} \left(\frac{1}{{\cal F}(s_{2,k})^2}-\frac{1}{{\cal
F}(s_{1,k})^2}\right) \ge \frac{1}{k_5}\sum_{k=1}^n (s_{2,k}-
s_{1,k}) \ge \frac{c}{k_5}s_{2,n}\,.
\]
This implies
\[
\frac{1}{{\cal F}(s_{2,n})^2}-\frac{1}{{\cal
F}(s_{1,1})^2}\ge\frac{c}{k_5}s_{2,n},
\]
since the intermediate terms are such that
\[
-\frac{1}{{\cal
F}(s_{1,k})^2}+\frac{1}{{\cal
F}(s_{2,k-1})^2}\le 0
\]
because $\mathcal{F}(s)$ is non-increasing and $s_{2,k-1}\le
s_{1,k}$. The monotonicity of the function ${\cal F}(s)$ allows then
to conclude that for all $s\in [s_{2,k}, s_{2,k+1}]$
\begin{equation}
\frac{1}{{\cal F}(s)^2}\ge \frac{1}{{\cal F}(s_{2,k})^2}\ge
\frac{1}{{\cal F}(s_{1,1})^2}+ \frac{c}{k_5}\,s_{2,k},
\end{equation}
Using the fact that $s\le s_{2,k+1}\le s_{2,k}+4$ we get
\begin{equation}
 {\cal F}(s)\le \frac{1}{\left[{\cal F}(s_{1,1})^{-2}+
c\,k_5^{-1}\,s\right]^{\frac{1}{2}}}\le
\frac{1}{(c_0+c_1s)^{\frac{1}{2}}}
\end{equation}
for large times $s$ and some positive constants $c_0, c_1$.
We have thus proved that the nonlinear entropy decays with the same
rate as the linear one,  when the initial relative mass is
nonzero.

\medskip

\noindent (ii) {\sc Proof without extra
restrictions.}
The arguments used  above and  in Appendix A4 are valid changing $h$ into $-h$
and $g$ into $-g$ under the same a priori bounds. Hence, the conclusion is valid
for negative and radial initial difference $v_0-V_{D_*}\le 0$.

To deal with the general case where $v_0-V_{D_*}$ is not radial or does
not have a sign, we use the maximum principle, after writing
$|v_0(x)-V_{D_*}(x)|\le  f(|x|)$. By comparison we have $v_1\le v\le
v_2$, where $v_1$ and $v_2$ are the solutions corresponding to
initial data $V_{D_*}- f$ and $V_{D_*}+ f$ resp. For the corresponding
$w=v/V_{D_*}$, $h=w-1$ and $g=h(D_{D_*}+y^2)$ a similar comparison holds.
Thus,  $w_1\le w\le w_2$, where $w_1\le 1$ and $w_2\ge 1$ are the
solutions with radial initial data $1 \pm (f/V_{D_*})$, hence functions
of $r=|y|$ and $s$. Same idea applies to $g$. Take now into account
the form of the entropy
\begin{equation}\label{Entropy.2}
{\mathcal F}[w]:= \frac1{1-m}\int_{\RR^d} \Psi(w)V_{{D_*}}^m\dy, \qquad
\mbox{with} \quad \Psi (w)=(w-1)-\frac 1m(w^m-1).
\end{equation}
We note that $\Psi(w)$ is convex and has a zero minimum  at $w=1$. Since we have just proved that
the decay result holds for both $g_1$ and $g_2$, the statement also
holds for $g$, even if we do not assume that $v_0-V_{D_*}$ is
nonnegative or radial.\qed

\medskip

%------------------------------------------------------------------------------
\noindent {\bf Proof of Corollary \ref{Conv.Weight}}. We recall
the following facts proved in \cite{BBDGV}, Lemma 6.2 under the
running assumptions, (H1) and (H2). First we have that for any $\vartheta\in[0,\frac{2-m}{1-m}]$,
there exists positive constants $K_\vartheta, K_2$ such that
\begin{equation*}
\left\||x|^\vartheta(v-V_{D_*})\right\|_{2}\leq K_\vartheta\left( \mathcal{F}[w]\right)^{1/2}\;.
\end{equation*}
Moreover
\begin{equation*}
\left\|v-V_{D_*}\right\|_{2}\leq K_2\left( \mathcal{F}[w]\right)^{1/2}\;.
\end{equation*}
We now recall the result of Lemma 3.6 of \cite{BBDGV}
\begin{equation}\label{Holder.Cont.Est}
\|v(s)-V_{D_*}\|_{C^\alpha(\RR^d)}\le\,\mathcal{H}\,\|v(s)-V_{D_*}\|_\infty\quad\forall\; t\geq t_0\;.
\end{equation}
for a suitable $\alpha\in(0,1)$, and we combine it with the interpolation inequality~\eqref{eq:interpolation}, with $\lambda=-\alpha\,d< 0=\mu<1/2=\nu$, $C=\mathcal C_{-\alpha d,\,0,\,1/2}$
\begin{equation*}
\|v(s)-V_{D_*}\|_\infty \le\,C\,\|v(s)-V_{D_*}\|_{C^\alpha}^\vartheta\; \|v(s)-V_{D_*}\|_2^{1-\vartheta} \le\,C\,\mathcal{H}^\vartheta\,\|v(s)-V_{D_*}\|_\infty^\vartheta\;\|v(s)-V_{D_*}\|_2^{1-\vartheta}
\end{equation*}
where $\vartheta=1/(2+\alpha\,d)$. This implies
\begin{equation*}
\|v(s)-V_{D_*}\|_\infty
    \le C^{1/(1-\vartheta)}\,\mathcal{H}^{\vartheta/(1-\vartheta)}\,\|v(s)-V_{D_*}\|_2
    \le \mathcal{K}_\vartheta \left( \mathcal{F}[w]\right)^{1/2}  \quad\forall\; t\geq t_0\;.
\end{equation*}
From H\"older's inequality,
\[
\|v(s)-V_{D_*}\|_q\le \|v(s)-V_{D_*}\|_\infty^{(q-2)/q}\, \|v(s)-V_{D_*}\|_2^{2/q}
    \le\mathcal{K}_q\left( \mathcal{F}[w]\right)^{1/2}
\]
for all $q\in[2,\infty]$, we deduce that $\|v(s)-V_{D_*}\|_q$ decays with the same rate as $\left( \mathcal{F}[w]\right)^{1/2}$.

If $q\in(1,2)$, we know from Lemma 6.2. of \cite{BBDGV} that there exists a positive constant $K(q)$ such that
\begin{equation*}
\|v-V_{D_*}\|_q\le K(q)\left( \mathcal{F}[w]\right)^{1/2},
\end{equation*}
This and the known decay of ${\mathcal F}$ proves (ii).

To prove (iii), use first \eqref{interp.Cj.Lp} with the choice $p=\infty$, i.e.
\begin{equation}\label{interp.Cj.Lp.infty}
\|f\|_{C^{j}(\RR^d)}\;\le\;\mathcal C_{j,d}\;\|f\|_{C^{j+1}(\RR^d)}^{\frac{j}{(j+1)}}\;\|f\|_\infty^{\frac{1}{j+1}}
\end{equation}
for any $j\in\mathbb{N}$, and the decay of the L$^\infty$ norm, namely $\|v-V_{D_*}\|_\infty\le Ks^{-1/4}$ to get
\begin{equation*}
\|v(s)-V_{D_*}\|_{C^j(\RR^d)}\le H_js^{-\frac{1}{4(j+1)}}\quad\forall\;s\geq s_0\;,
\end{equation*}
where in fact $H_j$ depends on $s$ itself and tends to zero as $s\to+\infty$, so that the bound can be improved.
Indeed we iterate the procedure putting such bound for the $C^{j+1}$ norm into \eqref{interp.Cj.Lp.infty} to get a new bound for the $C^{j}$ norm. In fact what we get after $h$ steps is, for any fixed $s\ge s_0$, $j\in{\mathbb N}$:
\[
\|v(s)-V_{D_*}\|_{C^j(\RR^d)}\le \frac{{\mathcal C}_{j,d}^{\sum_{0}^{h-1}\left(\frac j{j+1}\right)^h}k^{\frac 1j \sum_{1}^{h+1}\left(\frac j{j+1}\right)^h}H_j^{\left( \frac j{j+1}\right)^h}}{s^{k_h}}
\]
where the value of $k_h$ will be determined later. Notice in first place that the numerator of the above expression remain finite as $h\to\infty$, for any fixed $s\ge s_0$, $j\in{\mathbb N}$. As for $k_h$, by construction it satisfies the recursion relation
\[
k_0=\frac 1{4(j+1)},\ \ \ k_{h+1}=\frac j{j+1}k_h+\frac 1{4(j+1)}.
\]
subtracting $1/4$ to both sides of the latter equation gives
\[
k_{h+1}-\frac 14=\frac j{j+1}\left(k_h-\frac 14\right)
\]
which immediately gives $k_h=\frac 14-\left(\frac
j{j+1}\right)^h\frac j{4(j+1)}$, thus proving that $k_h\to 1/4$ as
$h\to+\infty$.\qed

\noindent{\bf Proof of Corollary \ref{thm:CRE-exp}.} We have
proved that $\mathcal{F}[w(s)]\le c_0s^{-1/2}$ and by Lemma
\ref{Lem.Bounds.RE} we also know that $F[w]\le c_1
\mathcal{F}[w]$. By Lemma \ref{Lem.Entr.Lp} with $m=m_*$, we have
\begin{equation}
\|w(s)-1\|_{d/2}^{d/2}=\|w(s)-1\|^{2+\frac{m_*}{1-m_*}}_{\LL^{2+\frac{m_*}{1-m_*}}(\RR^d)}
    \le \overline{D}_{m_*} F[w(s)] \le c_2 \mathcal{F}[w(s)] \le c_3s^{-1/2}
\end{equation}
Moreover, \eqref{infinity-interp} yields
\begin{equation}
\begin{split}
\|w(s)-1\|_{\LL^\infty(\RR^d)}\;
    &\le c_4\mathcal{F}[w(s)]^\frac{1}{d+2+\frac{m_*}{1-m_*}}=c_4\mathcal{F}[w(s)]^{2/(3d)}
     \le c_5 s^{-1/(3d)}\\
\end{split}
\end{equation}
Interpolating between these bounds shows that, for $q\in[d/2,+\infty]$:
\[
\|w(s)-1\|_q\le\|w-1\|_\infty^{[q-(d/2)]/q}\|w(s)-1\|_{d/2}^{d/(2q)}
        \le c_6s^{-\frac 13\left(\frac{1}{d}+\frac1{q}\right)}.
\]
To improve such bound we insert it in the interpolation inequality \eqref{interp.Cj.Lp} and use Theorem \ref{lem:holdereg} as well to get
\[
\|w(s)-1\|_{C^j(\RR^d)}\le \frac C{s^{\frac q3\left(\frac 1d+\frac 1q\right)\frac{k-j}{d+qk}}}
\]
for any $q\in[d/2,+\infty]$, $k>j\in{\mathbb N}$. As a function of $q$ the exponent of $s$ is nonincreasing, so we choose $q=d/2$ to get
\[
\|w(s)-1\|_{C^j(\RR^d)}\le \frac C{s^{\frac{k-j}{d(2+k)}}}.
\]
To optimize in $k$ we should take $k=\infty$, which is not allowed, so that for any fixed $\varepsilon>0$ we take $k$ large enough so that
\[
\|w(s)-1\|_{C^j(\RR^d)}\le \frac C{s^{\frac{1-\varepsilon}{d}}}
\]
as claimed. Putting this bound back into \eqref{interp.Cj.Lp} with $j=0$ and using what is known so far for the decay of the L$^p$ norm we get a decay of the form $\|w-1\|_\infty\le Cs^{-\alpha}$, $\alpha$ being given (with an inessential renaming of the free parameter $\varepsilon$) by
\[
\alpha(p,k)=\left(\frac{1-\varepsilon}{d}\right)\frac d{d+pk}+\frac13\left(\frac1d+\frac1p\right)\frac{pk}{d+pk}.
\]
Maximizing $\alpha$ w.r.t. to $p$ when $\varepsilon$ is small
enough yields again $p=d/2$, so that after some calculation we get
the exponent $\alpha(d/2,k)=\frac1d-\frac{2\varepsilon}{d(2+k)}.$ This
proves the claim for the L$^\infty$ norm and hence also for all
L$^q$ norms with $q\in(d/2,+\infty)$ by interpolation.\qed

%%%%%%%%%%%%%%%%%%%%%%%%%%%%%%%%%%%%%%%%%%%%%%%%%%%%%%%%%%%%%%%%%%%%%%%%%%

\section{Proofs for fast diffusion with $m\neq m_*$ revisited}\label{m.neq.mstar}

The previous method allows for shorter proofs of the convergence
when $m\ne m_*$, and at the same time some minor improvements of
paper \cite{BBDGV}. We recall that, in the case $m\not=m_*$, the
spectral gap inequality
\begin{equation}\label{HP.ineq}
F[g(s)]\le \lambda_{m,d}^{-1}\,I [g(s)]
\end{equation}
holds true, and the best constant is known for $m<m_*$, since \cite{BBDGV-CRAS, BBDGV}. We also recall the result of Proposition \ref{Fisher-lin-nonlin.2}
\[
I [g]=\int_{\RR^d}\left|\nabla g\right|^2 V_{D_*} \dy
    \le k_1\mathcal{I}[w]
      + k_3\mathcal{F}^{1+\sigma}[w]\\
\]
where, in particular, $k_1= 2W_1^{3-2m}$. From these bounds we get
\begin{equation}\label{E-Epr.m}
\begin{split}
\mathcal{F}(w)
    &\le^{(a)} \frac 12W_0^{m-2}\int_{\RR^d}|w-1|^2V_{D_*}^{m}\dy
        = \frac 12W_0^{m-2}F\\
    &\le^{(b)} \frac 12W_0^{m-2}\lambda_{m,d}^{-1}I [g] \\
    &\le^{(c)} \frac 12W_0^{m-2}\lambda_{m,d}^{-1}
        \left[k_1\mathcal{I}[w]
      + k_3\mathcal{F}^{1+\sigma}[w]\right]\\
\end{split}
\end{equation}
where in $(a)$ we compared the linear and nonlinear entropies via inequality \eqref{Disug.Entr.Lin-Nolin}, in $(b)$ the above spectral gap inequality, and in $(c)$ the above mentioned Proposition \ref{Fisher-lin-nonlin.2}. We may rewrite the latter formula as a differential inequality:
\[
\mathcal{F}' + W_0^{2-m}W_1^{2m-3}\lambda_{m,d}\mathcal{F}-\frac{W_1^{2m-3}}{2}k_3\mathcal{F}^{1+\sigma}\le 0
\]
so that by comparison with its explicit solution, we get
\[
\mathcal{F}[w(s)]
    \le \frac{\ee^{-k_4(s-s_0)}}
            {\left[\mathcal{F}[w(s_0)]^{-\sigma} + \frac{k_3 W_1^{2m-3}}{2k_4}\left(\ee^{-\sigma k_4 s}
                -\ee^{-\sigma k_4 s_0} \right) \right]^{\frac{1}{\sigma}}}
    \le k_5 \ee^{-k_4s}
\]
provided $\mathcal{F}[w(s_0)]$ is small enough, a property which holds for $s_0$ large enough.
This gives an exponential decay of the entropy, with
a rate $k_4=W_0^{2-m}W_1^{2m-3}\lambda_{m,d}$. Note that
$\lambda_{m,d}$ is the optimal constant in the Hardy- Poincar\'e
inequality \eqref{HP.ineq}, known for $m<m_*$ since \cite{BBDGV}.
Then we can proceed as in \cite{BBDGV} to show that the optimal
rate is given by $\lambda_{m,d}$. Indeed, one can substitute $W_0$
and $W_1$ with $\inf_{\RR^d}|w|$ and $\sup_{\RR^d}|w|$
respectively, allowing them to depend on time. Then prove that they
both tend to $1$ when $s\to +\infty$, so that
$k_4\to\lambda_{m,d}$; this can be done in view of the uniform
convergence of the relative error
$\|w(t)-1\|_{\LL^\infty(\RR^d)}\to 0$ as $s\to +\infty$ together
with a Gronwall-type argument. For more details we refer to
Section 6.3 of \cite{BBDGV}.

\medskip

%%%%%%%%%%%%%%%%%%%%%%%%%%%%%%%%%%%%%%%%%%%%%%%%
\noindent{\bf Remarks}

\noindent (i) When $m=m_*$ the above steps do not hold since we do
not have a spectral gap for the linearized generator. This is one of the reasons which forced
us to use Gagliardo-Nirenberg inequalities which, instead, compare
the Fisher information with a power of the entropy.

\noindent (ii) This method simplifies and complements some proofs
of \cite{BBDGV} when $m\neq m_*$, but also gives a more detailed
proof of the case $m\le 0$ that was only briefly treated in
\cite{BBDGV}\,. We finally emphasize the analysis of the present
paper covers the case $m=0$, that is logarithmic diffusion, even
in dimension $d=4$ since in that case $m_*=0$ so that no spectral
gap holds.

\noindent (iii) The interpolations made in the proof of Theorem
\ref{Conv.Weight} are valid also in the case $m\neq m_*$ and allow
to improve the convergence rate of the derivatives, proving that
the rate is always given by $\lambda_{m,d}$\,. We state here this
improved version of the main asymptotic Theorem of \cite{BBDGV}

\begin{thm}[Convergence with rate, $m\neq m_*$]
Under the assumptions of Theorem~{\rm~\ref{Thm:A1}}, if $m\ne m_*$, there exists $t_0\geq
0$ such that the following properties hold:
\begin{enumerate}
\item[{\rm (i)}] For any $q\in[q_*,\infty]$, there exists a
positive constant $C_q$ such that
\begin{equation*}
\|v(s)-V_{D_*}\|_q\le C_q\;\ee^{-\lambda_{m,d} \,s}\quad\forall\;s\geq s_0\;.
\end{equation*}
\item[{\rm (ii)}] For any $ \vartheta\in[0,(2-m)/(1-m)]$, there
exists a positive constant $K_\vartheta$ such that
\begin{equation*}
\big\|\,|x|^\vartheta(v(s)-V_{D_*})\big\|_{2}\le
K_\vartheta\;\ee^{-\lambda_{m,d}\,s}\quad\forall\;s\geq t_0\;.
\end{equation*}
\item[{\rm (iii)}] For any $j\in\mathbb{N}$, there exists a positive constant $H_j$ such that
\begin{equation*}
\|v(s)-V_{D_*}\|_{C^j(\RR^d)}\le H_j\,\ee^{-\lambda_{m,d}\,s}\quad\forall\;s\geq s_0\;.
\end{equation*}
\end{enumerate}
\end{thm}
%------------------------------------------------------------------------------
The constants $C_q$, $K_\vartheta$ and $H_j$ depend on $s_0$, $m$,
$d$, $v_0$, $D_0$, $D_1$, and $q$, $\vartheta$ and $j$; $s_0$ also
depends on $D_0$ and $D_1$. It is remarkable that the decay rate
of the nonlinear problem is given exactly by $\lambda_{m,d}$.
Rescaling back to the original equation, we obtain results in terms of intermediate
asymptotics, cf. Corollary \ref{Cor:A2} or Corollary 1.3 of \cite{BBDGV}\,.

\section*{Appendices}

\subsection*{\bf A1. Calculation of curvatures. Proof of Lemma
\ref{lemma.Rij}}

We can use well-known formulas for the Ricci tensor as a function
of the metric data:
$$
R_{ij}=g^{km}R_{ikjm}, \quad R_{ikjm}=\frac12\left(\partial^2_{kj}g_{im}+ \partial^2_{im}g_{kj}-
 \partial^2_{km}g_{ij}- \partial^2_{ij}g_{km} \right)+
 g_{np}(\Gamma^n_{kj}\Gamma^p_{im}+\Gamma^n_{km}\Gamma^p_{ij}).
$$
but in the case of conformal transformation there is a worked out relation between the
Ricci tensors of two metrics ${\bf g}$ and $\widetilde{{\bf g}}$ in terms of the conformal
factor relating them. Precisely, if   $\widetilde{{\bf g}}= (1/\varphi^2){\bf g}$, where $\varphi$ is a scalar,
the formula reads as follows \cite{Besse}:
$$
{\widetilde R}-R=\frac1{\varphi^2}\left[(d-2)\varphi\nabla^2\varphi+
(\varphi\Delta_{\bf g}\varphi - (d-1){\bf g}(\nabla \varphi, \nabla \varphi)
\cdot {\bf g}\right].
$$
where $R=(R_{ij})$ is the Ricci tensor of ${\bf g}$, ${\widetilde R}={\widetilde R}_{ij}$
is the Ricci tensor of ${\widetilde {\bf g}}$, $\nabla$ denotes
the gradient, $\nabla^2$ the Hessian and $\Delta_{\bf g}$ the
Laplace-Beltrami operator with respect to ${\bf g}$. Specializing the formula to the case
${\bf g}=\delta_{ij}$, so that $R_{ij}=0$, we get in coordinates
$$
{\widetilde R}_{ij}=
\frac1{\varphi^2}\left[(d-2)\varphi\partial^2_{ij}\varphi+
(\varphi\Delta\varphi - (d-1)|\nabla \varphi |^2)
\delta_{ij}\right].
$$
Put now ${\widetilde g}_{ij}=(1+|x|^2)^{-1}\delta_{ij}$ so that
$\varphi=(1+|x|^2)^{1/2}$. Then $\partial_i\varphi=x_i(1+|x|^2)^{-1/2},$
$$\partial^2_{ij}\varphi=-x_ix_j(1+|x|^2)^{-3/2}+ (1+|x|^2)^{-1/2}\delta_{ij},
\quad \Delta\varphi=-|x|^2(1+|x|^2)^{-3/2}+ d(1+|x|^2)^{-1/2}.
$$
Applying the last formula we get after
some calculations
%$$
%{\widetilde R}_{ij}=-(d-2)\frac{x_ix_j}{(1+|x|^2)^{2}}+
%\left[\frac{(d-2)(1+|x|^2)+ -|x|^2+ d(1+|x|^2)- (d-1)|x|^2}{(1+|x|^2)^{2}}\right]\delta_{ij}.
%$$
$$
{\widetilde R}_{ij}=-\frac{(d-2)x_ix_j}{(1+|x|^2)^{2}}+
\left[\frac{(d-2)|x|^2+ 2(d-1) }{(1+|x|^2)^{2}}\right]\delta_{ij}.
$$
There is a clear form of these expressions when we
take the particular point $\widehat{x}=(X,0,\cdot,0)$ which
implies no loss of geometrical generality since the metric is
conformal and radial, hence invariant under rotations in the
space. We get
\begin{equation}
\widetilde R_{11}(\widehat{x})=\frac{2(d-1)}{(1+X^2)^2}; \qquad
\widetilde R_{ii}(\widehat{x})=\frac{(d-2)X^2+2(d-1)}{(1+X^2)^2}
\quad \forall i=2,\cdots, d,
\end{equation}
and $\widetilde R_{ij}(\widehat{x})=0$ for all $i\ne j$. Both  eigenvalues
tend to zero as $|x|\to+\infty$ with different rates.  It immediately follows
that the symmetric tensor Ric is positive;
indeed, given $\xi\in {\mathbb R}^d$, we have
\[
\widetilde R_{ij}(\widehat{x})\xi_i\xi_j
\ge\frac{2(d-1)}{(1+X^2)^2}|\xi|^2>0,
\]
and the same is true for all $x\in \mathbb{R}^d$ by invariance
under rotations. If one wants to visualize the behaviour of this manifold, it is convenient to look
at Ricci curvatures given by
\begin{equation}
\widetilde r_1=\frac{\widetilde R(e_1,e_1)}{\widetilde g(e_1,e_1)}=\frac{2(d-1)}{(1+X^2)},
\qquad \widetilde r_i=\frac{\widetilde R(e_i,e_i)}{\widetilde g(e_i,e_i)}=\frac{2(d-1)+(d-2)X^2}{(1+X^2)},
\end{equation}
Note that the transversal curvatures tend to $(d-2)$ as $|x|\to \infty$ while the
curvature in the radial direction behaves like $O(|x|^{-2})$.
 This clearly shows the difference
in the behaviour of the curvatures in radial and transversal
directions which is typical of a cigar manifold.

Finally, the value of the scalar curvature follows from the
formula $R=g^{ij}R_{ij}.$ Since we are in a conformal situation,
it can be deduced in a direct way from the Yamabe formula
\cite{yamabe60}
$$
{\widetilde R}=-\frac{4(d-1)}{d-2}\frac{\Delta w}{w^{(d+2)/(d-2)}}, \quad
\mbox{with} \  w=g^{(d-2)/4},
$$
where ${\bf g}$ is the conformal factor, here $(1+|x|^2)^{-1}$,
cf. the formulas e.\,g. in \cite{VazSmooth}, pages 211-212.
In order to obtain the results stated in Lemma \ref{lemma.Rij} we only need to eliminate
the tildes from $\widetilde R_{ij}$ and $\widetilde R$. \qed

\medskip

\subsection*{\bf A2. Explicit representation of
the cigar}

We give here a simple parametric representation for the cigar-like
manifold $(M,{\bf g})$. We want to represent it as a hypersurface in $\RR^{d+1}$.
The radial symmetry of the metric suggests to represent such imbedded manifold
as $z=f(|y|)$ with variables $(y,z)\in \RR^{d+1}$, where $y\in\RR^d$ and $z\in\RR$, having a unique
chart $x\in \RR^d\mapsto (y,z)$  given by the formulas
\[
r=|y|=\Phi(\varrho)\qquad\mbox{and}\qquad z=\Psi(\varrho),
\]
where $\varrho=|x|\ge 0$. We fix $\Phi(0)=\Psi(0)=0$. We want the Euclidean metric in $\RR^{d+1}$
to induce the metric on the hypersurface. We know that the
infinitesimal length element in the radial direction satisfies
\[
\ds^2=\rd r^2 +\rd
z^2=\left[\Psi'^2(\varrho)+\Phi'^2(\varrho)\right]\rd\varrho^2
=\frac{\rd\varrho^2}{1+\varrho^2}
\]
which implies the relation $
\Psi'^2(\varrho)+\Phi'^2(\varrho)=1/(1+\varrho^2) $. On the other
hand, the length calculation for the transversal part gives
\[
\left[\frac{\Phi(\varrho)}{\varrho}\right]^2=\frac{1}{1+\varrho^2}.
\]
Solving the above two equations gives
\[
\Phi(\varrho)=\frac{\varrho}{(1+\varrho^2)^{1/2}}, \qquad
\Psi'(\varrho)=\frac{\varrho(2+\varrho^2)^{1/2}}{(1+\varrho^2)^{3/2}}.
\]
We see that  $r=\Phi(\varrho)$ goes from 0 to 1 as
$0<\varrho<\infty$. Analyzing the behaviour of $\Psi(\varrho)$,
one concludes that
\[
\Psi(\varrho)\approx \varrho^2 \quad\mbox{when }\varrho\approx 0,
\qquad \Psi(\varrho)\approx \log\varrho \quad\mbox{when
}\varrho\gg 1.
\]
This is the representation of a cigar.  The point out that the transversal radius
at infinity is  constant; actually, $\Phi(\varrho)\to 1$ as $\varrho\to +\infty$.

\medskip

\subsection*{\bf A3. Some Technicalities}

   We recall here some technical facts that we
used in the proofs. First we recall Theorem 2.4 of \cite{BBDGV},
\begin{thm}[Uniform $C^k$ regularity]\label{lem:holdereg}
Let $m <1$ and $w \in \LL^{\infty}_{{\rm loc}}((0,T) \times \RR^d)$ be a solution of~\eqref{eq.w}. Then for any $k\in\NN$, for any $s_0\in (0,T)$,
\begin{equation}
\sup_{s\geq s_0}\|w(s)\|_{C^k(\RR^d)}<+\infty\;.
\end{equation}
\end{thm}
We also needed an interpolation Lemma due to Gagliardo \cite{Ga}, cf. also  Nirenberg, \cite[p. 126]{MR0109940}.

\begin{lem} Let $\lambda$, $\mu$ and $\nu$ be such that $-\infty<\lambda \le \mu
\le \nu <\infty$. Then there exists a positive constant $\mathcal
C_{\lambda,\mu,\nu}$ independent of $f$ such that
\begin{equation}\label{eq:interpolation}
\|f\|_{1/\mu}^{\nu-\lambda} \le \mathcal
C_{\lambda,\mu,\nu}\|f\|_{1/\lambda}^{\nu-\mu} \;
\|f\|_{1/\nu}^{\mu-\lambda}\quad\forall\;f\in\mathcal C(\RR^d)\;,
\end{equation}
where $\|\cdot\|_{1/\sigma}$ stands for the following quantities:
\noindent(i) If $\sigma>0$, then
$\|f\|_{1/\sigma}=\left(\int_{\RR^d}|f|^{1/\sigma}\dx\right)^\sigma$.
\noindent(ii) If $\sigma<0$, let $k$ be the integer part of
$(-\sigma d)$ and $\alpha=|\sigma|d-k$ be the fractional
(positive) part of $\sigma$. Using the standard multi-index
notation, where $|\eta|=\eta_1+\ldots+\eta_d$ is the length
of the multi-index
$\eta=(\eta_1,\ldots\eta_d)\in\mathbb{Z}^d$, we define
\begin{equation*}\label{def.C^k}
\|f\|_{1/\sigma}=\left\{\begin{array}{lll} \displaystyle\max_{|\eta|=k}\; \big|\partial^\eta f\big|_\alpha=\displaystyle\max_{|\eta|=k}\; \sup_{x,y\in\RR^d}\;\dfrac{\big|\partial^\eta f(x)-\partial^\eta f(y)\big|}{|x-y|^\alpha}=|f\|_{C^\alpha(\RR^d)}& \mbox{if~}\alpha>0\;,\\[5mm]
\displaystyle\max_{|\eta|=k}\;\displaystyle\sup_{z\in\RR^d}\big|\partial^\eta f(z)\big|:=\|f\|_{C^k(\RR^d)}& \mbox{if~}\alpha=0\;.
\end{array}\right.
\end{equation*}
As a special case, we observe that $\|f\|_{-d/j}=\|f\|_{C^{j}(\RR^d)}$.

\noindent(iii) By convention, we note $\|f\|_{1/0}=\sup_{z\in\RR^d}|f(z)|=\|f\|_{C^0(\RR^d)}=\|f\|_{\infty}$.
\end{lem}

\medskip

\noindent\textbf{Remark.} The following special case
of the above interpolation inequality \eqref{eq:interpolation} has been used in the paper: let $k>j\in\mathbb{N}$ and $\lambda=-k/d\le\mu=-j/d\le\nu=1/p$.
Inequality~\eqref{eq:interpolation} becomes
\begin{equation}\label{interp.Cj.Lp}
\|f\|_{C^{j}(\RR^d)}\;\le\;\mathcal C_{j,k,p}\;\|f\|_{C^{k}(\RR^d)}^{\frac{d+jp}{d+kp}}\;\|f\|_p^{\frac{p(k-j)}{d+kp}}
\end{equation}
for any $k>j\in\mathbb{N}$ and $p>0$.

\medskip

\subsection*{\bf A4. Complete proof of the estimates of Section \ref{sect.nlem2}}

In the proof of Theorem \ref{thm.conv.entropy} in Section \ref{sect.nlem2}  we have assumed that
for every solution under the stated conditions there is an infinite sequence of intervals of {\sl good times} $[s_{1,k}\,,
s_{2,k}]\subset[2k,2k+2]$ with $s_{2,k}<s_{1,k+1}$ for all $k$, such
that
\begin{equation}\label{times.sk} I[g(s)]\le
k_4\,\mathcal{I}[g(s)]\qquad\mbox{\it for all } \ s\in
\bigcup_{k\in\NN}[s_{1,k}\,, s_{2,k}].
\end{equation}
Recall that in view of hypothesis (H2) and the discussion made in
Section \ref{sect.nlem2}, we may assume that $g$ is radially
symmetric and positive. We will also prove that the length, $l_k=
s_{2,k}- s_{1,k}$ of the intervals in our construction is at least
$1/2$ for all $k\ge k_0$.

The proof of these facts is long, and would have broken the flow of
the proof of Theorem \ref{thm.conv.entropy}: this is the reason why
we put it here. The main point in getting \eqref{times.sk} consists
in obtaining stronger estimates of the remainder term in the
inequality of Lemma \ref{Fisher-lin-nonlin} than the ones obtained
in Proposition \ref{Fisher-lin-nonlin.2}. We restate here Lemma 5.1
for convenience of the reader:

\textit{Let $0<W_0\le w\le W_1<+\infty$, be a measurable function on $\RR^d$,
with $W_0<1$ and $W_1>1$, and assume that $\mathcal{I}(w)<+\infty$. Then for
any $m<1$ the following inequality holds true
\begin{equation}
I [w(s)]
    \le k_1\mathcal{I}[w(s)]+ R[w(s)], \qquad \mbox{with } \ R[w(s)]=
      k_2\int_{\RR^d}g(s,y)^4V_*(y)^{4-3m}\dy\,,
\end{equation}
where $g=(w-1)V_*^{m-1}$; $k_1$ and $k_2$ are positive constants.}

We need to control the remainder term $R[w(s)]$ to proceed with the
asymptotic estimate. Note that $V_*^{4-3m}$ is integrable for
$m=m^*$: in fact, for such a value of $m$ we have
\begin{equation}\label{Fisher.Reminder}
R[w(s)]= k_2\int_{\RR^d}\frac{g(s,y)^4}{(1+y^2)^{(d+4)/2}}\dy
\end{equation}
Put now $N(s)=N[g(s)]=\|g(\cdot,s)\|_\infty$, the supremum of $g$
for fixed time $s>0$. We know that $N(s)$ is  uniformly bounded in
time. Then we have
\begin{equation}\label{better}
 R[w(s)]\le k_3 \,[N(s)]^4
\end{equation}
We want to estimate the decay of $N(s)$ in time in terms of the
linearized Fisher information $I [w]$. Suppose for a moment that  we
can prove that the remainder term is small relatively to $I [w]$,
more precisely that $R[w(s)]\le \frac12I [w(s)]$ for all large $s$.
In that case we conclude that $I [w(s)] \le 2 k_1\mathcal{I}[w(s)],$
and the desired estimate \eqref{times.sk} easily follows. Hence, we
need to prove that
\begin{equation}\label{conditionK}
N(s)^4 \le \frac{I[g(s)]}{K}
\end{equation}
with $K>k_3$, say $K=2k_3$. This is a most convenient estimate on
the values of $g$.

Unfortunately, even under such assumption it is not clear that the
last inequality holds at all times, or even at all large times.
Therefore, we shall be cautious and call {\it good times } those
times at which \eqref{conditionK}  holds with $K\ge 2k_3$. The
frequency and density of the intervals of such times is important,
as the end of proof of Theorem \ref{thm.conv.entropy} shows.

We will now proceed with the proof of the existence of the time
intervals stated at the beginning of this section. They will consist
only of so-called good times. The proof is split into two parts,
namely

\noindent(i) \textit{Controlling the remainder term away from the
origin.} This is the part where we use the fact the
$|v_0-V_*|\le \tilde f$ for a {\it radially symmetric} $\tilde f$;

\noindent(ii) \textit{Transforming the outer control into a control
on a small ball,} namely we will control $\sup_{r\ge R} g(r,s)$ for
small $R>0$. Due to the peculiarities of the parabolic Harnack
inequality, we shall prove that such a control only takes place for
a large set of so-called good times, this sufficing for our goals.

\noindent \textbf{Part (i). The control of a radial $g$ far away from the origin}\\
In the calculation that follows we drop the $s$-dependence for
convenience since time does not enter in the argument. Let $g(r)$ be
a \textit{nonnegative} continuous function such that $g(r)\to 0$ as
$r\to\infty$ and let
\[
{\cal M}_R=\int_R^{\infty}\frac{g(r)}{(1+r^2)^{d/2}}r^{d-1}\rd
r<\infty\,,\qquad\mbox{and}\qquad
I_R=\int_R^{\infty}\frac{|g'(r)|^2}{(1+r^2)^{(d-2)/2}}r^{d-1}\rd
r<\infty\,.
\]
We put the powers in a way such that it is clear that for $r>1$ we are dealing with the
radial case and we are merely asking the mass and the linearized Fisher information to be finite.
Now pick $\alpha>0$, $R_1>R>1$ and calculate
\[\begin{split}
g^{1+\alpha}(R)-g^{1+\alpha}(R_1)
    &=-\alpha\int_R^{R_1}g'(r)g^\alpha(r)\rd r
     \le \alpha\int_R^{R_1}|g'(r)|g(r)^\alpha\rd r
     \le \alpha\int_R^{R_1}|g'(r)|r^{\frac{1}{2}}\frac{g(r)^\alpha}{r^{\frac{1}{2}}}\rd r\\
    &\le \alpha\left[\int_R^{R_1}|g'(r)|^2 r\rd r\right]^{\frac{1}{2}}
        \left[\int_R^{R_1}\frac{g(r)^{2\alpha}}{r}\rd r\right]^{\frac{1}{2}}
\end{split}
\]
Now, if we assume $\alpha\ge 1/2$, the last integral can be bounded as follows:
\[
\int_R^{R_1}\frac{|g(r)|^{2\alpha}}{r}\rd r
    \le \sup_{r\ge R}|g(r)|^{2\alpha-1}\,
        \int_R^{\infty}\frac{|g(r)|}{r}\rd r,
\]
so that for $\alpha\ge 1/2$ we have obtained
\[\begin{split}
g^{1+\alpha}(R)-g^{1+\alpha}(R_1)\le \alpha \left[\sup_{r\ge 1}|g(r)|^{2\alpha-1}\,
        \int_1^{\infty}\frac{|g(r)|}{r}\rd r \right]^{\frac{1}{2}}
    \left[\int_1^{\infty}|g'(r)|^2 r\rd r\right]^{\frac{1}{2}}\,.
\end{split}
\]
Letting $R_1\to\infty$, and assuming that $g(R)\to 0$ as $R\to \infty$, we get
\[\begin{split}
g^{1+\alpha}(R)\le \alpha \left[\sup_{r\ge 1}|g(r)|^{2\alpha-1}\,
        \int_1^{\infty}\frac{|g(r)|}{r}\rd r \right]^{\frac{1}{2}}
    \left[\int_1^{\infty}|g'(r)|^2 r\rd r\right]^{\frac{1}{2}}\,.
\end{split}
\]
Taking the supremum over $R\ge1$ on the l.h.s. and simplifying we get:
\[\begin{split}
\left[\sup_{R\ge 1}g(R)\right]^{4}\le \alpha^{\frac{8}{3}}
    \left[\int_1^{\infty}\frac{|g(r)|}{r}\rd r \right]^{\frac{4}{3}}
    \left[\int_1^{\infty}|g'(r)|^2 r\rd r\right]^{\frac{4}{3}}\le
    c\,{\cal M}_1^{\frac{4}{3}}\,I_1^{\frac{4}{3}}\le\tilde c I^{\frac43}.
\end{split}
\]
This is a very good estimate because it says that the supremum of
$g^4$ outside the unit ball is not only proportional to $I$ as
expected in so-called better times, but even more: it is proportional to a
higher power of $I$. Now, recall that $I[w(s)]\to 0$ as $s\to\infty$.
If the same could be done near $r=0$ the proof that every large $s$
is a good time would be complete.

The previous calculation can be done in the complement of the ball of
radius $R$ as small as we like and then $g(R)$ will depend also on
an inverse power of $R$, because of the presence of the factors $1+r^2$
in the denominators of the last quantities.
We now get in the last line for $0<R<1$
\[\begin{split}
\left[\sup_{r\ge R}|g(r)|\right]^{4}\le \alpha^{\frac{8}{3}}
    \left[\int_R^{\infty}\frac{|g(r)|}{r}\rd r \right]^{\frac{4}{3}}
    \left[\int_R^{\infty}|g'(r)|^2 r\rd r\right]^{\frac{4}{3}}\le
    \frac{C}{R^{8(d-1)/3}} {\cal M}_R^{\frac{4}{3}}\,I_R^{\frac{4}{3}}.
\end{split}
\]
The estimate blows up at $R=0$. Therefore, we cannot let $R\to 0$ to
get an estimate for  $N(s)$ for any $x\in \RR^d$.

\noindent\textit{Justifying that $g$ goes to zero at infinity.} To
conclude part (i) of the proof, it remains to prove that $g(R)\to 0$
as $R\to \infty$. Choose $R_n$ such that
\[
\int_{R_n}^{\infty}|g'|^2 r\,\rd r<\frac{1}{4n^2}\qquad\mbox{and}\qquad
\int_{R_n}^{\infty}\frac{|g|^2}{r} \,\rd r<\frac{1}{4n^2}
\]
and define $\widetilde R_n=\min\left\{r\ge R_n\;; g^2(r)\le \frac{1}{2n^2}\right\}$.
Indeed the set in the r.h.s. is not empty since $g/r$ is integrable at infinity: this is not compatible with $g$ being everywhere larger than a positive constant for all $r\ge R_n$. Notice that $\widetilde R_n\ge R_n$ and $g^2(\widetilde R_n)\le1/(2n^2)$. Hence, for all $R\ge R_n$:
\[
g^2(R)=g^2(\widetilde R_n)+2\int_{\widetilde R_n}^{R}g\,g'\rd r \le \frac{1}{2n^2} +2\left|\int_{\widetilde R_n}^{R}g|g'|\rd r\right|
        \le \frac{1}{2n^2} +
2\left[\int_{R_n}^{\infty}|g'|^2\rd r\right]^{\frac{1}{2}}\left[\int_{R_n}^{\infty}g^2\rd r\right]^{\frac{1}{2}}\le \frac1{n^2}
\]
Therefore, $0\le g(R)\le 1/n$ for all $R\ge R_n$. The proof of part (i) is now complete.

\medskip

\noindent\textbf{Part (ii). Transforming the outer control into a
control on a small ball.}

\noindent In part (i) we have
estimated the supremum of a radial $g^4$ outside a ball of radius
$R>0$ in terms of $I[g(s)]^{4/3}$, so that the problem is to
estimate the supremum  inside a ball as well, hopefully in terms of
$I[g(s)]^{1+\alpha}$, at least in the form $\ve I[g(s)]$. We are
unable to prove that for all (sufficiently large) times. To
circumvent such a difficulty we have to make use of a rather
complicated argument that takes into account the possibility that
such estimate does not hold because of possible bad behaviour of $g$
at points near the origin. We begin by carefully labeling the times. We say that a time
$s\in[0,\infty)$ belongs to the class of {\sl good times} $\mathcal{G}_K$, if
\[
N(s)^4=\sup_{y}(g(s,y))^4 < \frac{I[g(s)]}{K}
\]
We are not claiming that some half-line
$[T,\infty)\subset {\cal G}_{2k_3}$, which would finish the proof in
the simple way. Finally, we say that a time is {\sl very good},
$s\in {\cal V}_C$, if
\begin{equation}
N(s)\le C \, I[g(s)]^{4/3}
\end{equation}
for some $C>0$, in the spirit of the radial estimate away from the
origin. Note that since $I(s)\to 0$ we have the inclusion of very
good times with constant $C$ into the good times with any constant
$K>0$ if $s$ is large enough.

\medskip

\noindent\textsc{Harnack inequality.} The study of points near the
space origin is based on classical regularity theory for linear or
quasilinear parabolic equations in divergence form. We are going to use the version of
the celebrated paper by Aronson-Serrin \cite{AS} . We consider
the equation satisfied by the error function
\begin{equation}
h(s,y)=w(s,y)-1=\frac{v(s,y)}{V_{D_*}(y)}-1.
\end{equation}
It can be written in the standard form
\begin{equation}
\partial_sh=\nabla\cdot {\bf A}(y,h,\nabla h) + B(y,h,\nabla h).
\end{equation}
In fact, starting with the equation satisfied by $w$, we have
\begin{equation}\label{eq.w-1}
\partial_s h  = \partial_s w     =\frac{1}{V_*}\nabla\cdot\left[w  V_*\nabla
        \left(\frac{w ^{m-1}-1}{m-1}V_*^{m-1}\right)\right]
    =\frac{1}{V_*}\nabla\cdot\left[(h +1) V_*\nabla
        \left(\frac{(h +1)^{m-1}-1}{m-1}V_*^{m-1}\right)\right]
\end{equation}
so that we can identify
\[
{\bf A}=(h+1)\nabla\left(\frac{(h+1)^{m-1}-1}{m-1}V_*^{m-1}\right),\
\ \ B=(h+1)\frac{\nabla V_*}{V_*}
\cdot\nabla\left[\frac{(h+1)^{m-1}-1}{m-1}V_*^{m-1}\right]
\]
and we have to check that the structure conditions are satisfied by
${\bf A}$ and $B$  in a compact ball with constants that do not depend on $s$.
In fact, the structure conditions are satisfied in the homogeneous form of \cite{AS},
which means that, in the notation of that paper, the terms $f,g,h=0$ in the structure
condition for ${\bf A}, B$ vanish. Checking this is a straightforward
calculation involving also the known bounds on $w$ namely $0<W_0\le
w=h+1\le W_1$. We note in passing that since we already know that  $w\to 1$ uniformly in $\RR^d$ as $s\to\infty$,
 the lower and upper bounds  $W_0,W_1$ can be taken closer and closer to 1 if we restrict the time to
 $s\ge s_0$ and $s_0$ is large enough.
%\[\begin{split}
%|B| &=\left|V_*^{m-2}(h+1)^{m-1}\nabla h\cdot\nabla V_* + \left[(h+1)^{m-1}-1\right](h+1)V_*^{m-2}|\nabla V_*|^2\right|\\
%    &\le |\nabla V_*|V_*^{m-2}(h+1)^{m-1}|\nabla h|+ \underbrace{4V_*(h+1)\left[(h+1)^{m-1}-1\right]}_{\le c h}
%\end{split}\]

In any case, we conclude that the Harnack inequality has the standard form,
as stated below. This also implies a similar Harnack inequality for
$g$ if we work on a bounded space domain, say, in $B_1(0)$. We state it next. Take $T>0$ large and consider the parabolic cylinders
\[
Q=[T-2,T]\times B_1(0)\,,\qquad Q_{1/2}=[T-1/2, T]\times B_{1}(0)\,,\qquad \widetilde{Q}=[T-2, T-1]\times B_1(0)\,.
\]
The parabolic Harnack inequality on the disjoint cylinders $Q_{1/2}$
and $\widetilde{Q}$ is then
\begin{equation}
\inf_{Q_{1/2}}g(s,y)\ge c\; \sup_{\widetilde{Q}}g(s,y)=c\;
\widetilde{N}
\end{equation}
for some positive constant $c<1$ depending only on structural constants.
Note that since $g$ is continuous on $Q$ all the above suprema are attained at some points.

\medskip

%%%%%%%%%%%%%%%%%%%%%%%%%%%%%%
\noindent\textsc{Evolution of the maximum of $h$.} The equation for
$h$ can be written as
\[\begin{aligned}
\partial_sh&=\nabla\cdot\left[(h+1)\nabla\left(\frac{(h+1)^{m-1}-1}{m-1}V_*^{m-1}\right)\right]+V_*^{m-2}(h+1)^{m-1}\nabla V_*\cdot\nabla h\\
&+V_*^{m-3}|\nabla V_*|^2(h+1)[(h+1)^{m-1}-1]=\\
&=\nabla h\cdot\nabla\left(\frac{(h+1)^{m-1}-1}{m-1}V_*^{m-1}\right)+(h+1)\nabla\cdot[(h+1)^{m-2}V_*^{m-1}\nabla h]\\
&+mV_*^{m-2}(h+1)^{m-1}\nabla V_*\cdot\nabla h\\
&+(h+1){[(h+1)^{m-1}-1]}{\nabla\cdot(V_*^{m-2}\nabla V_*)}
 +{V_*^{m-3}|\nabla V_*|^2}(h+1)[(h+1)^{m-1}-1]
\end{aligned}
\]
In particular, using the fact $h$ is small we get
\begin{equation}
\partial_sh= \mbox{second and first order terms in $h$}+  C(r,h)\,h
\end{equation}
where $C(r,h)\le \overline{k}$ for a suitable $\overline{k}$
independent of $r$ and depending only on the known a priori bounds
for $h$. Then, as a consequence of the Maximum Principle, cf. e.g. \cite{AS}, the
maximum $N(s)$ obeys the growth rate
\begin{equation}
N'(s)\le \overline{k} N(s).
\end{equation}
In fact the function $H(s,y)=N(s_0)e^{\overline{k}(s-s_0)}$ is an
explicit supersolution in the whole space for $s\ge s_0$.

\medskip

\noindent\textsc{The structure of good times. Alternative.} Now we
are going to prove that

\begin{lem} For every time interval $(T-2,T)$ with $T$ large enough there
is at least a subinterval of length 1/2 consisting of very good
times.\end{lem}

 \noindent {\sl Proof.~} We can use the cylinders $Q_{1/2}$ and $\widetilde{Q}$ as in
the previous paragraph on the Harnack inequality. The idea is to consider separately the two
possibilities: (a) either the maximum in $x$ of $h$ at every time of
the lower cylinder is taken outside the ball $B_2(0)$, or (b) the
maximum at one of such times, say $s_0$, is taken inside.

In case \textsl{(b)}, $T$ must  be a good time if it is large
enough as we show next. We take $s_0\in (T-2,T-1)$ as above and let
$N_1$ be the corresponding maximum in the ball. Of course $N_1\le
\widetilde N$. We now apply the growth rate of previous
paragraph to obtain
\begin{equation}
N(s_2)\le C\,N_1\le C\widetilde{N}, \qquad C=e^{2\overline{k}}.
\end{equation}
for every $s_2$ in the upper cylinder: $T-(1/2)\le s_2\le T$. On the
other hand, for every $y\in B_1(0)$ we have for such $s_2$ the lower
estimate $h(s_2,y)\ge c\widetilde{N}$. We conclude that the maximum
and the minimum at all those times are related by a constant. This
is also true for the function $g$ up to a small change in the
constant, hence
\begin{equation}
g(s_2,y'_M)\le C_1 g(s_2,y)
\end{equation}
where now $y_M'$ is the point of maximum of $g$ in the ball
$B_1(0)$. Now, we know that on the boundary $|y|=1$ there is a good
estimate for $g$, more precisely, for such $y$ of unit norm  $g(s_2,y)$ satisfies the
estimates that defines the {\sl very good time}, and it does with a fixed constant $C$.
We conclude that $g(s_2,y)$, $|y|\ge 1$, also satisfies such an estimate with a possible worse constant
$C'=C_1C$. Since the estimate was true for $|y|\ge 1$ we are done.

 Therefore, whenever $T$ is  not a good  time in ${\cal G}_K$ with $K=2k_3$, the first part of
the alternative, \textsl{(a)}, must be true. But in that case for every time in the
lower cylinder we know that the maximum of $h$ is taken in the
exterior of the ball $B_1$. If we look for the expression of
$g=h\,(1+|y|^2)$  this also means that the maximum of $g(s,\cdot)$ at
times in $(T-2,T-1)$ is  taken at an exterior point (maybe
different). So all these times are very good. Recalling what was said before,
they are good times in ${\cal G}_K$ if $T$ is large enough.\qed

 \noindent\textsc{Choice of intervals of
good times.} We can apply the previous results letting $T=2k+2$, for
$k\ge k_0$ and $k_0$ sufficiently large. The above lemma implies that there
exists a subinterval $[s_{1,k}\,, s_{2,k}]\subset[2k,2k+2]$ of
length at least $1/2$ made of times in ${\cal G}_K$ with $K=2k_3$.

%%%%%%%%%%%%%%%%%%%%%%%%%%%%%%%%%%%%%%%%%%%%%%%%%%%%%%%%
\section{Concluding remarks and open problems}
\label{sec.cr}

{\bf\ref{sec.cr}.1.} The special situation has been studied for the
critical exponent $m_*=(d-4)/(d-2)$ in dimensions $d\ge 3$ where our
considerations make sense. Algebraic extensions have been shown to
be fruitful or intriguing in some dynamical studies. Here, for $d=2$
we formally have $m_*= \pm\infty$, which is an extreme situation for
porous medium that has appeared in the literature (for instance, in
connection with the mesa problem), \cite{mesa, CF, VazSmooth}, while
for $d=1$ we formally get $m_*=3$, a value inside the porous medium
range where nothing special has been shown to happen.

\medskip

\noindent {\bf \ref{sec.cr}.2.} We pose the following questions:

$\bullet$ Are the rates obtained in this paper optimal for a
certain class of data, as the linearized analysis suggests?

$\bullet$ Can we prove convergence,  maybe with worse rates or
without rates, for more general initial data? we recall that for
$m>m_c$ all nonnegative initial data in $L^1(\RR^d)$ are attracted
towards a Barenblatt solution, with no rate in that generality.

  $\bullet$ Find an explicit optimal dependence of the
constant in the asymptotic formula with respect to the data.

$\bullet$ Assuming that we get an optimal rate of convergence, can
we find a profile for the next level of approximation?

\noindent {\bf\ref{sec.cr}.3.} One may wonder if the techniques used in \cite{BBDGV-CRAS,BBDGV} for the case $m\neq m_*$, which use Hardy-Poincar\'e inequalities, work also in the case $m_*$. We have partially given a negative answer to this question in Corollary \ref{No.Hardy.No.Party}, in which we have shown that no inequality of Hardy type can hold for the linearized Fisher information $I[w]$. However, one may wonder if modified
versions of the Hardy-Poincar\'e inequalities, with logarithmic terms added in the spirit of the classical
Hardy inequality in $\RR^2$, allow to solve the problem. Thus, there is a family of valid Hardy inequalities (see below) in which the Dirichlet form involved has a logarithmic correction, but we are not able to prove the
asymptotic results by means of such inequalities. It is then a further open problem to see whether this path may lead to the goal or not.

\begin{prop}\label{HP.log}
Let $d\ge 3$. We have
\begin{equation}\label{HP.log.ineq}
\int_{\RR^d}g^2\,\rd\mu_\alpha\,\le
\mathcal{H}_{\alpha\,,d}\,\int_{\RR^d}\left|\nabla
g\right|^2\,\rd\nu_\alpha
\end{equation}
for any $g\in\mathcal D(\RR^d)$ and for any
$0<\alpha\le\frac{d}{2}-1$, where
\begin{equation}\begin{split}
\rd\mu_\alpha(y)&=\left(1+|y|^2\right)^{-\frac{d}{2}}\,\left[1+\log(1+|y|^2)\right]^{\alpha-1}\,\rd y,\\
\rd\nu_\alpha(y)&=\left(1+|y|^2\right)^{1-\frac{d}{2}}\,\left[1+\log(1+|y|^2)\right]^{\alpha+1}\rd y\\
\end{split}
\end{equation}
and
\begin{equation}\label{Const.Gap.Lem}
\mathcal{H}_{\alpha\,,d}=\frac{2(d-2)}{\alpha\,(d-2-2\alpha)\,\min
\Big\{2\alpha\,,(d-2-2\alpha)\Big\}}.%=\left\{\begin{array}{lll}
%\frac{(d-2)}{\alpha^2\,(d-2-2\alpha)} &\mbox{if}~~0<\alpha\le\frac{d-2}{4}
%\\[3mm]
%\frac{2(d-2)}{\alpha\,(d-2-2\alpha)^2}~~&\mbox{if}~~\alpha>\frac{d-2}{4} \\
%\end{array}
%\right.
\end{equation}
\end{prop}

\noindent{\bf Acknowledgements.} G.G. acknowledges an ESF grant of
the research group GLOBAL (Global and Geometric Aspects of Partial
Differential equations), which allowed him to visit M.B. and J.L.V.
at the Universidad Aut{\'o}noma de Madrid, and the Departamento de
Matem\'{a}ticas of such University, where part of this work has been
done.  M.B. and J.L.V. were partially supported by Spanish Project
MTM2005-08760-C02-01. M.B. and J.L.V. wish to thank the Dipartimento
di Matematica of Politecnico di Torino where part of the work has
been done.  We thank J. Dolbeault for a careful reading
of this article which resulted in a serious improvement.

\

%%%%%%%%%%%%%%%%%%%%%%%%

\end{document}